\documentclass[a4paper,11pt]{article}
\usepackage{amsmath,latexsym,amssymb,amsfonts,natbib}
\usepackage[a4paper,hmargin=1.87 cm]{geometry}
\usepackage[dvips]{graphicx}
\usepackage{pifont}
\usepackage{pstcol,pst-fill,pst-grad}
\usepackage{rotating}
\usepackage{subfigure}
\usepackage{graphics}
\usepackage{hyperref}
\usepackage{multirow}
\usepackage{bigints}
\usepackage{multicol}
\usepackage{courier}
\usepackage[latin1]{inputenc}
\usepackage{indentfirst}
\usepackage{ccaption}
\usepackage{verbatim}
\usepackage{makeidx}
\usepackage{listings}
\usepackage{dsfont}
\usepackage[toc,page]{appendix}

\newcommand{\proof}     {\paragraph{Proof}}
\newcommand{\carre}     {\hfill$\Box$}
\numberwithin{equation}{section}
\newtheorem{defi}{Definition}[section]
\newtheorem{lem}{Lemma}[section]

\newtheorem{prop}{Proposition}[section]
\newtheorem{rem}{Remark}

\newtheorem{ass}{Assumption}

\newtheorem{ex}{Example}[section]

\title{Local stability of McKean-Vlasov equations arising from heterogeneous Gibbs systems using limit of relative entropies}
\author{ Donald A. Dawson\thanks{School of Mathematics and Statistics, Carleton University, 1125 Colonel By Drive
Ottawa, Ontario K1S 5B6, Canada. ddawson@math.carleton.ca}, Ahmed Sid-Ali\thanks{School of Mathematics and Statistics,
Carleton University, 1125 Colonel By Drive Ottawa, Ontario K1S 5B6, Canada. ahmedsidali@cunet.carleton.ca}, Yiqiang Q. Zhao\thanks{School of Mathematics and Statistics, Carleton University, 1125 Colonel By Drive Ottawa, Ontario K1S 5B6, Canada. zhao@math.carleton.ca} \\
{\small }\\
{\small }\\
}

\date{}
\begin{document}
\maketitle

\begin{abstract}
A family of heterogeneous mean-field systems with jumps is analyzed. These systems are constructed as a Gibbs measure on block graphs. When the total number of particles goes to infinity, a law of large numbers is shown to hold in a multi-class context resulting in the weak convergence of the empirical vector towards the solution of a McKean-Vlasov system of equations. We then investigate the local stability of the limiting McKean-Vlasov system through the construction of a local Lyapunov function. We first compute the limit of adequately scaled relative entropy functions associated with the explicit stationary distribution of the $N$-particles system. Using a Laplace principle for empirical vectors we show that the limit takes an explicit form. Then we demonstrate that this limit satisfies a descent property which, combined with some mild assumptions shows that it is indeed a local Lyapunov function.   
\end{abstract}


\section*{Introduction}
    
The study of heterogeneous mean-field systems is a growing area of research. The main motivation is that the homogeneity assumption underlying the classical mean-field models often no longer holds when considering systems outside statistical physics. Therefore, many researchers studied systems with different heterogeneous assumptions; see, e.g., \cite{Contucci+al2008, Daw+Sid+zha2020, Kno+Lowe+Schub2021, Lack+Ram2019} and the references therein for an overview of the recent development in the subject. The focus in the current paper is on a particular family of mean-field systems of Gibbs type constructed on block graphs. This family is a particular instance of the models introduced in \cite{Daw+Sid+zha2020, Daw+Sid+zha2021} with the particularity laying in the specification of a heterogeneous Gibbs measure as a stationary distribution. This was in particular inspired by the interacting particle systems of Gibbs type on complete graphs analyzed in \cite{Budh+Du(a)2015, Budh+Du(b)2015}.

Classical questions in the study of mean-field systems include their asymptotic behavior when the total number of particles $N$ in the system and/or the time $t$ tends to infinity. Moreover, under which conditions one can justify the interchange of the limits $N\rightarrow\infty$ and $t\rightarrow\infty$. In particular, we will show for the specific family of Gibbs systems introduced in Section \ref{Gibb-sect} that, under mild assumptions, the associated empirical vector converges weakly and uniformly over compact time intervals, as $N\rightarrow\infty$, towards the solution of a McKean-Vlasov system of equations (see Theorem \ref{LLN-prop}). This kind of result is known in the literature as \textit{law of large numbers}. Thus, as a consequence, the McKean-Vlasov limiting system can be used to approximate the behavior of the large particles system over finite time intervals. Thence, one might wonder whether or not this approximation is still relevant when time goes to infinity. This question turns out to be related to the stability properties of the limiting McKean-Vlasov system. More precisely, if the latter contains a unique asymptotically stable equilibrium, then the interchange $N\rightarrow\infty$ and $t\rightarrow\infty$ is fully justified. However, if the McKean-Vlasov system contains multiple equilibria or a unique equilibrium that is not stable, care must be taken since the approximation is not necessarily accurate. The intuition behind this is that if the McKean-Vlasov limit system has several $\omega$-limit sets, the question is which of these sets characterizes the large-time behavior of the system. In addition, in such a case, metastable phenomena might be observed resulting in transitions of the empirical vector process from one $\omega$-limit set to another. For a more detailed discussion, one can consult \cite{Daw+Sid+zha2021, Ben+Leboud2008} and the references therein. Therefore, the stability of the limiting system is of great interest, which motivated us to investigate it in the current paper. 

The classical approach for studying the stability of dynamical systems is through the construction of Lyapunov functions. However, for McKean-Vlasov equations, given the nonlinearity, finding Lyapunov functions in the general case is very challenging. Nevertheless, the Gibbs nature of the systems studied here allows us to adopt the limit of relative entropy approach introduced in \cite{Budh+Du(a)2015, Budh+Du(b)2015}. The main idea behind this approach takes root in the observation that the relative entropy is a natural Lyapunov function for linear ergodic Markov processes; see, e.g. \cite{Spitzer71}. Moreover, despite the nonlinearity of the McKean-Vlasov system, the corresponding $N$-particle processes describe a linear Markov process. Therefore, one then proposes the limit of suitably normalized relative entropies associated with the stationary distribution of the $N$-particles system as the candidate Lyapunov function, providing that the limit takes a useful form. This approach was shown in \cite{Budh+Du(a)2015, Budh+Du(b)2015} to be successful for a family of homogeneous mean-field models with jumps of Gibbs type. The goal of the current paper is to extend the approach to the multi-class setting tackled here.

Notice that for the general non-Gibbs family of mean-field models on block graphs introduced in \cite{Daw+Sid+zha2020, Daw+Sid+zha2021}, the stationary distribution usually will not take an explicit form, and thus a different approach is needed. One possible line of work is the approach proposed in \cite{Budh+Du(a)2015} for particular systems with exchangeable particles, where instead of considering limits of relative entropies associated with the stationary distribution, one may consider the large time $t$ and large particles $N$ limits associated with the exchangeable joint probability distribution. Therefore, given the multi-class structure of the systems in \cite{Daw+Sid+zha2020, Daw+Sid+zha2021}, one might consider specific systems with multi-exchangeable particles. Another possible approach to construct Lyapunov functions is through the {\it Friedlin-Wentzell quasipotential} \cite{Freid+Wentz2012}. Indeed, a close tie between the quasipotential associated with small noise stochastic systems and Lyapunov functions for the underlying deterministic models was observed in the literature; see, e.g., \cite{DeMarco+Qian1989, DeMarco+Berg1987}. Therefore, for the general models introduced in \cite{Daw+Sid+zha2020, Daw+Sid+zha2021}, one might view the finite $N$-particles system as a small noise perturbation of the limiting McKean-Vlasov system. Notice that the specific quasipotential for these models was introduced in \cite{Daw+Sid+zha2021}. The idea is then to investigate under which condition the quasipotential can serve as a Lyapunov function for the McKean-Vlasov system. This however goes beyond the scope of the current paper. 

The rest of the paper is organized as follows. We introduce in Section \ref{Gibb-sect} the family of Gibbs systems on block graphs and show some preliminary properties. We then prove in Section \ref{LLN-sect} the law of large numbers and the convergence of the $N$-particles empirical vector toward the solution of a McKean-Vlasov system of equations. Therefore, we investigate in Section \ref{stab-sect} the local stability of the limiting McKean-Vlasov system by constructing a candidate Lyapunov function.  Using a particular Laplace principle associated with the vector of empirical measures (see Proposition \ref{Sanov-vect-meas-prop}), we start by computing in Proposition \ref{limit-rel-entr-prop} the limit of suitably normalized relative entropies and show that the limit takes an explicit form. Then, we show in Lemma \ref{crit-point-lem} that the limiting function characterizes the fixed points of the McKean-Vlasov system. Finally, Proposition \ref{descent-prop-lyap} shows that this limiting function satisfies a descent property, which, combined with mild assumptions, shows that it is indeed a local Lyapunov function for the McKean-Vlasov system of equations.

\section{Gibbs measures on block graphs}
\label{Gibb-sect}
The Gibbs measure concept has a long history and plays an important role in statistical physics. The underlying principle is that, when a system is in equilibrium, states with lower energy are more likely than those with higher energy. Thence, J.W. Gibbs proposed the probability measure $exp\{- U ({\pmb x}) / \kappa T \} $ to capture this principle, where $ \kappa $ is the Boltzmann constant, $T$ is the temperature and the function $U({\pmb x})$ gives the energy of the system when it is in the state ${\pmb x}$. Thence, the Gibbs measure being a model of equilibrium, given an energy function, one might seek a Markov process for which the Gibbs measure is the stationary distribution. Such Markov processes are often referred to as {\it Glauber dynamics} thanks to his seminal paper \cite{Glauber63}. For a detailed introduction to the topic, one can consult, e.g., \cite{Marti99, Stroock2005}.  

 Consider a graph $\mathcal{G}=(\mathcal{V},\Xi)$ composed of $r$ blocks $C_1,\ldots,C_r$ of sizes, respectively, $N_1,\ldots,N_r$, where $\mathcal{V}$ is the set of the nodes and $\Xi$ is the set of the edges. Denote by $|\mathcal{V}|=N_1+\cdots+N_r=N$ the total number of nodes in the graph. Suppose that each block $C_j$ is a clique, i.e., all the $N_j$ nodes of the same block are connected. Furthermore, the nodes within the same block $C_j$ are decomposed into two subsets: 
      \begin{itemize}
       \item \textbf{The set of central nodes $C^c_j$:} composed of the nodes that are connected to all the other nodes within the same block but not to any node from the other blocks. We set $|C_j^c|=N_j^c$.
       \item \textbf{The set of peripheral nodes $C^p_j$:}  composed of the nodes that are connected to all the other nodes within the same block and to all the peripheral nodes from the other blocks. We set $|C^p_j|=N^p_j$. Thus, the sub-graph engendered by all the peripheral nodes in $\mathcal{G}$ is complete.  
      \end{itemize}  
      
 The graph $\mathcal{G}=(\mathcal{V},\Xi)$ is thus composed of $2r$ components. This will play an important role in the upcoming analysis. We associate each node of the graph with a particle taking values in a finite state space $\mathcal{Z}=\{1,2,\ldots,K\}\subset\mathbb{N}$. Denoting by ${\pmb x}=(x_n,x_m,n\in C_j^c,m\in C_j^p,1\leq j \leq r)\in\mathcal{Z}^N$ a configuration of the $N$ particles over the graph  $\mathcal{G}=(\mathcal{V},\Xi)$, the corresponding local empirical measures describing the state of each component are defined by
\begin{align}
\varrho^{j,\iota,N}_z({\pmb x})=\frac{1}{N_j^{\iota}}\sum_{i\in C_j^{\iota}}\mathds{1}_{\{x_i=z\}},\quad \mbox{for any $z\in\mathcal{Z}$, $\iota\in\{c,p\}$, and $1\leq j\leq r$}.
\label{loc-emp-meas-gibbs}
\end{align}

Notice that for each $1\leq j\leq r$ and $\iota\in\{c,p\}$, the local empirical measure $\varrho^{j,\iota,N}_z({\pmb x})$ takes values in the state spaces $\mathcal{M}^{N_j^{\iota}}_1(\mathcal{Z})=\mathcal{M}_1(\mathcal{Z})\cap\frac{1}{N_j^{\iota}}\mathbb{Z}^K$, where $\mathcal{M}_1(\mathcal{Z})$ is the space of probability measures on $\mathcal{Z}$ endowed with the topology of weak convergence and $\mathbb{Z}$ is the set of integers. The corresponding empirical vector describing the state of the entire system is denoted by 
\begin{align*}
\varrho^N({\pmb x})=(\varrho^{1,c,N}({\pmb x}),\varrho^{1,p,N}({\pmb x}),\ldots,\varrho^{r,c,N}({\pmb x}),\varrho^{r,p,N}({\pmb x}))\in\prod_{j=1}^r\prod_{\iota\in\{c,p\}}\mathcal{M}^{N_j^{\iota}}_1(\mathcal{Z})\subset (\mathcal{M}_1(\mathcal{Z}))^{2r}.
\end{align*}
For ease of readability, we suppress in the sequel the dependency of $\varrho^{j,\iota,N}$ and $\varrho^N$ upon $N$. Thus we simply write $\varrho^{j,\iota}$ and $\varrho$ instead. We associate with the configuration ${\pmb x}$ of the $N$-particles the following \textit{energy function} 
 \begin{align}
U_N({\pmb x})=\sum_{i=1}^N V(x_i)+ \frac{\beta}{2N} \sum_{j=1}^r\bigg[\sum_{i\in C_j^c}\sum_{k\in C_j}W(x_i,x_k)+\sum_{i\in C_j^p}\sum_{k\in C_j^c\cup C^p} W(x_i,x_k)\bigg],
\label{energy-func}
\end{align}
 where $V:\mathcal{Z}\rightarrow\mathbb{R}$ is the \textit{potential function}, $W:\mathcal{Z}\times \mathcal{Z}\rightarrow\mathbb{R}$ is the symmetric \textit{interaction function}, and $\beta>0$ is the \textit{interaction parameter}. Thence, the corresponding Gibbs measure is given by 
\begin{align}
\pi^N({\pmb x})=\frac{1}{Z_N}\exp\{-U_N({\pmb x})\},
\label{stat-distr-gibbs}
\end{align}
where ${Z_N}$ is a normalization constant. We will now construct a Markov process, namely a {\it Glauber dynamics}, having the Gibbs measure $\pi^N$ as its stationary distribution. To this end, let us first introduce the directed graph $(\mathcal{Z},\mathcal{E})$ with $\mathcal{E}\subset\mathcal{Z}\times\mathcal{Z}\backslash \{(z,z)| z \in\mathcal{Z}\}$ representing the set of admissible jumps. Thence, whenever $(z,z')\in\mathcal{E}$, a particle at state $z$ is allowed to transit to $z'\in\mathcal{Z}$ at a rate that depends on the current state of the node and the state of its neighbors. Before going further, suppose the following assumptions throughout the paper.
\begin{ass} 
\label{ass-prin}
\begin{enumerate}
\item The set of edges $\mathcal{E}$ is symmetric. 
\item The directed graph $(\mathcal{Z},\mathcal{E})$ is irreducible.
\item For each block $1\leq j\leq r$, there exist $p_j^c,p_j^p,\alpha_j\in(0,1)$ such that, as $N\rightarrow\infty$,
\begin{align}
\frac{N_j}{N}\rightarrow \alpha_j,\quad \frac{N_j^p}{N_j}\rightarrow p_j^p,\quad \frac{N_j^c}{N_j}\rightarrow p_j^c,\quad p^p_j+p^c_j=1,\quad\text{and}\quad\sum_j\alpha_j=1.
\label{p-regul}
\end{align}
\end{enumerate}
\end{ass}
 
Define the matrix $(\alpha(z,z'))_{z,z'\in\mathcal{Z}}$ that identifies the allowed transitions for one particle as 
 \begin{align*}
 \alpha(z,z')=\left\{\begin{tabular}{c c} 
                1&\mbox{if $(z,z')\in\mathcal{E}$,}\\
                0&\mbox{else}.
 \end{tabular}\right.
 \end{align*}
 Moreover, let the matrix $A_N({\pmb x},{\pmb y})$ indexed by the elements ${\pmb x},{\pmb y}\in\mathcal{Z}^N$ be defined as
 \begin{align*}
 A_N({\pmb x},{\pmb y})=\left\{\begin{tabular}{c c} 
                $a(x_l,y_l)$&\mbox{if ${\pmb x}$ and ${\pmb y}$ differ exactly in one index $l\in C_j$, for some $1\leq j\leq r$}.\\
                0&\mbox{else}, 
 \end{tabular}\right.
 \end{align*}
Hence $A_N$ determines which states of the $N$-particle system can be reached in one jump. At this stage, there are several ways to construct Glauber dynamics. See, e.g., \cite[Sect. 3]{Marti99} for an overview. We propose here the {\it Metropolis dynamic} characterized by the following rate matrix: 
 \begin{align}
 \psi^N({\pmb x}, {\pmb y})= e^{-\big(U^N ({\pmb y})-U^N ({\pmb x})\big)^+} A_N ({\pmb x},{\pmb y})\quad \mbox{for ${\pmb x},{\pmb y} \in\mathcal{Z}^N,{\pmb x}\neq {\pmb y}$},
\label{metro-dyn-mat} 
 \end{align}
and $ \psi^N({\pmb x}, {\pmb x})=-\sum_{{\pmb y}\neq {\pmb x}} \psi^N({\pmb x}, {\pmb y})$, for all ${\pmb x} \in\mathcal{Z}^N$. One can verify that the rate matrix $\psi^N$ has $\pi^N$ as its stationary distribution. In fact, consider two configurations ${\pmb x},{\pmb y} \in\mathcal{Z}^N$. By symmetry, one has $A_N ({\pmb x}, {\pmb y}) = A_N ({\pmb y}, {\pmb x})$. Moreover, using $(\ref{stat-distr-gibbs})$ and $(\ref{metro-dyn-mat})$, it is easily to check that $\pi^N ({\pmb x})\psi^N({\pmb x}, {\pmb y})= \pi^N ({\pmb y})\psi^N({\pmb y}, {\pmb x})$. Then $\psi^N({\pmb x}, {\pmb y})$ satisfies the detailed balance condition with respect to $\pi^N$. Hence, $\pi^N$ is a stationary distribution of the Markov chain, and the rate matrix $\psi^N$ is reversible with respect to $\pi^N$. Furthermore, since the graph of allowed transitions $(\mathcal{Z},\mathcal{E})$ is irreducible by Assumption $\ref{ass-prin}$, the rate matrix $\psi^N$ is also irreducible and thus, $\pi^N$ is the unique stationary distribution. 
  
One might observe from the rate matrix introduced in $(\ref{metro-dyn-mat})$ that the transition between two configurations ${\pmb x}$ and ${\pmb y}$ depends on the difference between their total energy. Therefore, to investigate the large-scale behavior of the system, we first estimate this difference when the total number $N$ of particles in the system is very large. Given the multi-class structure of the system, we take throughout the paper the convention that as $N\rightarrow\infty$, $\min\limits_{\substack{1\leq j\leq r \\ \iota\in\{c,p\}}}N_j^{\iota}\rightarrow\infty$. In addition, and for the aim of simplicity, we will ignore from now on the environment potential by supposing that $V\equiv 0$ and thus focus on the interaction component of the system. Nevertheless, our results can be easily extended to the case with non-zero potential. 

Let us define, for $x,y\in\mathcal{Z}$, $q=(q^{1,c},q^{1,p},\ldots,q^{r,c},q^{r,p})\in (\mathcal{M}_1(\mathcal{Z}))^{2r}$, and $a,b_1,\ldots,b_r\in\mathbb{R}$, the following real-valued functions 
\begin{align}
\psi^{j,c,N}( x,y,q,a,b_1)&=\frac{\beta}{N}\bigg[a\sum_{z\in\mathcal{Z}}\big(W(z,y)-W(z,x)\big) q^{j,c}_z+b_1\sum_{z\in\mathcal{Z}}\big(W(y,z)-W(x,z)\big)q^{j,p}_z\bigg],
\label{psi-c-func}
\end{align}
\begin{align}
\psi^{j,p,N}( x,y,q,a,b_1,\ldots,b_r)&=\frac{\beta}{N}\bigg[a\sum_{z\in\mathcal{Z}} \big( W(y,z)-W(x,z)\big)q^{j,c}_z
+\sum_{l=1}^r\bigg(b_l\sum_{z\in\mathcal{Z}} \big( W(z,y)-W(z,x)\big)q^{l,p}_z\bigg)\bigg],
\label{psi-p-func}
\end{align} 
\begin{align*}
B^{c,N}(x,y)=\frac{\beta}{2N}\bigg(W(y,y)+W(x,x)-2W(x,y)\bigg),\quad\text{and}\quad B^{p,N}(x,y)=\frac{\beta}{2N}\bigg(W(x,x)- W(y,x)\bigg).
\end{align*} 

 \begin{lem}
 Let ${\pmb x},{\pmb y} \in\mathcal{Z}^N$ be two configurations such that $A_N({\pmb x},{\pmb y})=1$. If the unique index satisfying $x_l \neq y_l$ is a central node, i.e., $l\in C_{j^*}^c$, for some block $1\leq j^*\leq r$, then
\begin{align}
U_N({\pmb y})-U_N({\pmb x})=\psi^{j^*,c,N}(x_l,y_l,\varrho({\pmb x}),N^c_{j^*},N^p_{j^*})+B^{c,N}(x_l,y_l),
\label{ener-func-diff-c}
\end{align}
 for which there exists a constant $C\in (0,\infty)$ such that 
  \begin{align}
\sup_{x,y\in\mathcal{Z}}|B^{c,N}(x,y)| \leq\frac{C}{N}.
\label{B-func-c-bound}
  \end{align}
  
If the unique index such that $x_l \neq y_l$ is a peripheral node, i.e. $l\in C_{j^*}^p$, for some block $1\leq j^*\leq r$, then  
 \begin{align}
 U_N({\pmb y})-U_N({\pmb x})=\psi^{j^*,p,N}(x_l,y_l,\varrho({\pmb x}), N^c_{j^*},N_1^p,\ldots, N_r^p)+B^{p,N}(x_l,y_l),
 \label{ener-func-diff-p}
 \end{align}   
for which there exists a constant $C'\in (0,\infty)$ such that 
\begin{align}
\sup_{x,y\in\mathcal{Z}}|B^{p,N}(x,y)| \leq\frac{C'}{N}.
\label{B-func-p-bound}
\end{align}
\label{lem-glaub-dyn}
 \end{lem}

\proof Let ${\pmb x},{\pmb y} \in\mathcal{Z}^N$ be two configurations such that $A_N({\pmb x},{\pmb y})=1$, and let $l$ be the unique index such that $x_l\neq y_l$. The index $l$ can either refer to a central or a peripheral node. We treat the two cases separately. 
 
\paragraph{Case 1.} Suppose that the two vectors ${\pmb x}$ and ${\pmb y}$ differ in one central node, i.e., $l\in C_{j^*}^c$ for some $1\leq j^*\leq r$. Therefore by $(\ref{energy-func})$ one obtains
\begin{equation}
\begin{split}
U_N({\pmb y})-U_N({\pmb x})&=\frac{\beta}{2N} \bigg(\sum_{i\in C_{j^*}^c}\sum_{k\in C_{j^*}}\big(W(y_i,y_k)-W(x_i,x_k)\big)+\sum_{i\in C_{j^*}^p}\sum_{k\in C_{j^*}^c\cup C^p} \big( W(y_i,y_k)-W(x_i,x_k)\big)\bigg).
\end{split}
\end{equation}
By the symmetry of the interaction function $W$ and using (\ref{loc-emp-meas-gibbs}) one further obtains
\begin{equation}
\begin{split}
\sum_{i\in C_{j^*}^c}\sum_{k\in C_{j^*}}\big(W(y_i,y_k)-W(x_i,x_k)\big)&=\sum_{i\in C_{j^*}^c}\bigg[\sum_{k\in C^c_{j^*}}\big(W(y_i,y_k)-W(x_i,x_k)\big)+\sum_{k\in C^p_{j^*}}\big(W(y_i,y_k)-W(x_i,x_k)\big)\bigg]\\
&=\sum_{i\in C_{j^*}^c, i\neq l}\bigg[\sum_{k\in C^c_{j^*}}\big(W(y_i,y_k)-W(x_i,x_k)\big)+\sum_{k\in C^p_{j^*}}\big(W(y_i,y_k)-W(x_i,x_k)\big)\bigg]\\
&\qquad+\bigg[\sum_{k\in C^c_{j^*}}\big(W(y_l,y_k)-W(x_l,x_k)\big)+\sum_{k\in C^p_{j^*}}\big(W(y_l,y_k)-W(x_l,x_k)\big)\bigg]\\
&=\sum_{i\in C_{j^*}^c, i\neq l}\bigg[\sum_{k\in C^c_{j^*},k\neq l}\big(W(y_i,y_k)-W(x_i,x_k)\big)+\big(W(y_i,y_l)-W(x_i,x_l)\big)\bigg]\\
&\qquad+\sum_{k\in C^c_{j^*},k\neq l}\big(W(y_l,y_k)-W(x_l,x_k)\big)+\big(W(y_l,y_l)-W(x_l,x_l)\big)\\
&\qquad+\sum_{k\in C^p_{j^*}}\big(W(y_l,y_k)-W(x_l,x_k)\big)\\
&=2N_{j^*}^c\sum_{z\in\mathcal{Z}}\big(W(z,y_l)-W(z,x_l)\big)\varrho^{j^*,c}_z({\pmb x})-2 W(x_l,y_l)+W(x_l,x_l)+W(y_l,y_l)\\
&\qquad+N_{j^*}^p\sum_{z\in\mathcal{Z}}\big(W(y_l,z)-W(x_l,z)\big)\varrho^{j^*,p}_z({\pmb x}),\\
\end{split}
\end{equation}
and
\begin{equation}
\begin{split}
\sum_{i\in C_{j^*}^p}\sum_{k\in C_{j^*}^c\cup C^p} \big( W(y_i,y_k)-W(x_i,x_k)\big)&=\sum_{i\in C_{j^*}^p}\bigg[\sum_{k\in C_{j^*}^c} \big( W(y_i,y_k)-W(x_i,x_k)\big)+\sum_{k\in  C^p} \big( W(y_i,y_k)-W(x_i,x_k)\big)\bigg]\\
&=\sum_{i\in C_{j^*}^p}\bigg[\sum_{k\in C_{j^*}^c,k\neq l} \big( W(y_i,y_k)-W(x_i,x_k)\big)+\big( W(y_i,y_l)-W(x_i,x_l)\big)\bigg]\\
&=N_{j^*}^p\sum_{z\in \mathcal{Z}}\big( W(z,y_l)-W(z,x_l)\big)\varrho^{{j^*},p}_z({\pmb x}).
\end{split}
\end{equation}
Therefore, one concludes that
\begin{equation}
\begin{split}
U_N({\pmb y})-U_N({\pmb x})&=\frac{\beta}{2N}\bigg[2N_{j^*}^c\sum_{z\in\mathcal{Z}}\big(W(z,y_l)-W(z,x_l)\big)\varrho^{j^*,c}_z({\pmb x})+2N_{j^*}^p\sum_{z\in\mathcal{Z}}\big(W(y_l,z)-W(x_l,z)\big)\varrho^{j^*,p}_z({\pmb x})\\
&\qquad\qquad\qquad+W(y_l,y_l)+W(x_l,x_l)-2W(x_l,y_l)\bigg],\\
\end{split}
\end{equation}
from which $(\ref{ener-func-diff-c})$ and $(\ref{B-func-c-bound})$ follow.

\paragraph{Case 2:} Suppose now that the two vectors ${\pmb x}$ and ${\pmb y}$ differs in one peripheral node, i.e., $l\in C_{j^*}^p$ for some $1\leq j^*\leq r$. Therefore from $(\ref{energy-func})$ one gets
 \begin{equation}
\begin{split}
U_N({\pmb y})-U_N({\pmb x})&=\frac{\beta}{2N} \bigg[\sum_{i\in C_{j^*}^c}\sum_{k\in C_{j^*}}\big(W(y_i,y_k)-W(x_i,x_k)\big)+\sum_{i\in C_{j^*}^p}\sum_{k\in C_{j^*}^c\cup C^p} \big( W(y_i,y_k)-W(x_i,x_k)\big)\bigg]\\
&\qquad+\frac{\beta}{N} \sum_{j=1,j\neq j^*}^r\bigg[\sum_{i\in C_{j}^p}\sum_{k\in C_{j}^c\cup C^p} \big( W(y_i,y_k)-W(x_i,x_k)\big)\bigg].
\end{split}
\end{equation}
Again, using the symmetry of $W$ together with (\ref{loc-emp-meas-gibbs}) one further obtains
\begin{equation}
\begin{split}
\sum_{i\in C_{j^*}^c}\sum_{k\in C_{j^*}}\big(W(y_i,y_k)-W(x_i,x_k)\big)&=N^c_{j^*}\sum_{z\in\mathcal{Z}}\big(W(z,y_l)-W(z,x_l)\big)\varrho^{j^*,c}_z({\pmb x}),
\end{split}
\end{equation}
\begin{equation}
\begin{split}
\sum_{i\in C_{j^*}^p}\sum_{k\in C_{j^*}^c\cup C^p} \big( W(y_i,y_k)-W(x_i,x_k)\big)&=2N^p_{j^*}\sum_{z\in\mathcal{Z}}\big( W(z,y_l)-W(z,x_l)\big)\varrho^{{j^*},p}_z({\pmb x})\\
&\qquad+N^c_{j^*}\sum_{z\in\mathcal{Z}} \big( W(y_l,z)-W(x_l,z)\big)\varrho^{{j^*},c}_z({\pmb x})\\
&\qquad+\sum_{j=1,j\neq j^*}^r\bigg[N^p_{j}\sum_{z\in\mathcal{Z}} \big( W(y_l,z)-W(x_l,z)\big)\varrho^{j,p}_z({\pmb x})\bigg]\\
&\qquad+\big( W(y_l,y_l)-W(x_l,x_l)\big),
\end{split}
\end{equation}
and
\begin{equation}
\begin{split}
\sum_{j=1,j\neq j^*}^r\bigg[\sum_{i\in C_{j}^p}\sum_{k\in C_{j}^c\cup C^p} \big( W(y_i,y_k)-W(x_i,x_k)\bigg]=\sum_{j=1,j\neq j^*}^r\bigg[N_j^p\sum_{z\in\mathcal{Z}} \big( W(z,y_l)-W(z,x_l)\varrho^{j,p}_z({\pmb x})\bigg].
\end{split}
\end{equation}
Thus, one deduces that
\begin{align*}
U_N({\pmb y})-U_N({\pmb x})&=\frac{\beta}{2N}\bigg[2N^c_{j^*}\sum_{z\in\mathcal{Z}} \big( W(y_l,z)-W(x_l,z)\big)\varrho^{{j^*},c}_z({\pmb x})+2\sum_{j=1}^r\bigg(N_j^p\sum_{z\in\mathcal{Z}} \big( W(z,y_l)-W(z,x_l)\varrho^{j,p}_z({\pmb x})\bigg)\\
&\qquad\qquad\qquad +\big( W(y_l,y_l)-W(x_l,x_l)\big)\bigg],
\end{align*}
from which $(\ref{ener-func-diff-p})$ and $(\ref{B-func-p-bound})$ follow. 
\carre 
\\
\\ 
One might notice from Lemma \ref{lem-glaub-dyn} that if two configurations differ in only one index $l$, then the jump rates of the Markov process governed by the rate matrix $\psi^N$ depend on the states $x_i$ of the other particles $i\neq l$ only through the local empirical measures $\varrho^{1,c}({\pmb x}),\varrho^{1,p}({\pmb x}),\ldots,\varrho^{r,c}({\pmb x}),\varrho^{r,p}({\pmb x})$. To emphasize this fact, we introduce the rate functions $\lambda^{j,c,N}$ and $\lambda^{j,p,N}$ defined, for all $1\leq j\leq r$, $z,z'\in\mathcal{Z}$, $q=(q^{1,c},q^{1,p},\ldots,q^{r,c},q^{r,p})\in (\mathcal{M}_1(\mathcal{Z}))^{2r}$, and $a,b_1,\ldots,b_r\in\mathbb{R}$, by 
 \begin{align}
\lambda^{j,c,N}_{z,z'}\big(q,a,b_1\big)= e^{-\big(\psi^{j,c,N}(z,z',q,a,b_1)+B^{c,N}(z,z')\big)^+} a(z,z'),
 \label{rate-funct-cent-N} 
 \end{align}
 \begin{align}
\lambda^{j,p,N}_{z,z'}\big(q,a,b_1,\ldots,b_r\big)= e^{-\big(\psi^{j,p,N}(z,z',q,a,b_1,\ldots,b_r)+B^{p,N}(z,z')\big)^+} a (z,z').
\label{rate-funct-periph-N} 
 \end{align} 
Thus, when the system is in configuration ${\pmb x}$, a central node $l \in C_{j}^c$ of a given block $1\leq j\leq r$ jumps from a state $z$ to $z'$ at rate 
 \begin{align}
\lambda^{j,c,N}_{z,z'}\big(\varrho({\pmb x}),N^c_{j},N^p_{j}\big),
\label{rate-funct-cent-N-2} 
 \end{align}
and a peripheral node $l \in C_{j}^p$ within the block $1\leq j\leq r$ jumps from a state $z$ to $ z'$ at rate 
 \begin{align}
\lambda^{j,p,N}_{z,z'}\big(\varrho({\pmb x}),N^c_{j},N_1^p,\ldots, N_r^p\big).
\label{rate-funct-periph-N-2} 
 \end{align}

Now, using $(\ref{psi-c-func})$ and $(\ref{psi-p-func})$ together with Assumption \ref{ass-prin}, one can easily prove that the functions $\psi^{j,c,N}(\cdot,N^c_{j},N^p_{j})$ and $\psi^{j,p,N}(\cdot,N^c_{j},N_1^p,\ldots, N_r^p)$ converge respectively, as $N\rightarrow\infty$, toward the functions $\psi^{c}$ and $\psi^{p}$ defined, for all $1\leq j\leq r$, $z,z'\in\mathcal{Z}$ and $q=(q^{1,c},q^{1,p},\ldots,q^{r,c},q^{r,p})\in (\mathcal{M}_1(\mathcal{Z}))^{2r}$ by
\begin{align*}
\psi^{j,c}( z,z',q)&=\beta\bigg[\alpha_j p_j^c\sum_{x\in\mathcal{Z}}\big(W(x,z')-W(x,z)\big) q^{j,c}_x+\alpha_j p_j^p\sum_{x\in\mathcal{Z}}\big(W(z',x)-W(z,x)\big)q^{j,p}_x\bigg],\\
\psi^{j,p}( z,z',q)&=\beta\bigg[\alpha_j p_j^c\sum_{x\in\mathcal{Z}} \big( W(z',x)-W(z,x)\big)q^{j,c}_x
+\sum_{l=1}^r\bigg(\alpha_l p_l^p\sum_{x\in\mathcal{Z}} \big( W(x,z')-W(x,z)\big)q^{l,p}_x\bigg)\bigg].
\end{align*}
Using this together with the inequalities in $(\ref{B-func-c-bound})$ and $(\ref{B-func-p-bound})$, one can further prove that the rate functions $\lambda^{j,c,N}_{z,z'}\big(\cdot,N^c_{j},N^p_{j}\big)$ and $\lambda^{j,p,N}_{z,z'}\big(\cdot,N^c_{j},N_1^p,\ldots, N_r^p\big)$ converge respectively, as $N\rightarrow\infty$, toward the functions $\lambda^{j,c}_{z,z'}$ and $\lambda^{j,p}_{z,z'}$ defined by
 \begin{align}
\lambda^{j,c}_{z,z'}(q)= e^{-\big(\psi^{j,c}(z,z',q)\big)^+} a(z,z')\mbox{ and } \lambda^{j,p}_{z,z'}(q)&= e^{-\big(\psi^{j,p}(z,z',q)\big)^+} a (z,z').
\label{lambd-gibbs}
 \end{align}

Fix $q=(q^{1,c},q^{1,p},\ldots,q^{r,c},q^{r,p})\in (\mathcal{M}_1(\mathcal{Z}))^{2r}$ then, for all $1\leq j\leq r$, the rate matrix $A^{j,c}_q:=\big(\lambda^{j,c}_{z,z'}(q)\big)_{(z,z')\in\mathcal{Z}\times\mathcal{Z}}$, with $\lambda^{j,c}_{z,z}(q)=-\sum_{z'\neq z}\lambda^{j,c}_{z,z'}(q)$, is the generator of an ergodic Markov chain with a unique invariant measure defined, for all $z\in\mathcal{Z}$, by 
\begin{align}
\pi^{j,c}(q)_z=\frac{1}{Z_N(q)}\exp\bigg\{-\beta \big[\alpha_jp_j^c\sum\limits_{x\in\mathcal{Z}}W(x,z)q_x^{j,c}+\alpha_jp_j^p\sum\limits_{x\in\mathcal{Z}}W(x,z)q_x^{j,p}\big]\bigg\},
\label{stat-dist-fix-c}
\end{align}
where $Z_N(q)$ is a normalization constant given by 
\begin{align*}
Z_N(q)=\sum_{z\in\mathcal{Z}}\exp\big\{-\beta \big[\alpha_jp_j^c\sum\limits_{x\in\mathcal{Z}}W(x,z)q_x^{j,c}+\alpha_jp_j^p\sum\limits_{x\in\mathcal{Z}}W(x,z)q_x^{j,p}\big]\big\}.
\end{align*}
In fact, this can be easily verified by checking the detailed balance condition. Similarly, the rate matrix $A^{j,p}_q:=\big(\lambda^{j,p}_{z,z'}(q)\big)_{(z,z')\in\mathcal{Z}\times\mathcal{Z}}$, with $\lambda^p_{z,z}(q)=-\sum_{z'\neq z}\lambda^{j,p}_{z,z'}(q)$, is the generator of an ergodic Markov chain with a unique stationary distribution defined, for all $z\in\mathcal{Z}$, by 
\begin{align}
\pi^{j,p}(q)_z=\frac{1}{Z_N(q)}\exp\bigg\{-\beta \big[\alpha_jp_j^c\sum\limits_{x\in\mathcal{Z}}W(x,z)q_x^{j,c}+\sum\limits_{l=1}^r\alpha_lp_l^p\sum\limits_{x\in\mathcal{Z}}W(x,z)q_x^{l,p}\big]\bigg\},
\label{stat-dist-fix-p}
\end{align}
where again $Z_N(q)$ is the normalization constant given by 
\begin{align*}
Z_N(q)=\sum_{z\in\mathcal{Z}}\exp\bigg\{-\beta \big[\alpha_jp_j^c\sum\limits_{x\in\mathcal{Z}}W(x,z)q_x^{j,c}+\sum\limits_{l=1}^r\alpha_lp_l^p\sum\limits_{x\in\mathcal{Z}}W(x,z)q_x^{l,p}\big]\bigg\}.
\end{align*}

\section{Law of large numbers and McKean-Vlasov limiting system}
\label{LLN-sect}
 Denote by $\big({\pmb x}^N(t)\big)_{t \geq 0}=\big(x_n(t),x_m(t),n\in C_j^c, m\in C_j^p,1\leq j\leq r\big)_{t\geq 0}$ the stochastic process characterizing the state of the system at each time $t\geq 0$. Hence $\big({\pmb x}^N(t)\big)_{t \geq 0}$ is a Markov process with state space $\mathcal{Z}^{2r}$ and transition rate matrix $\big(\psi^N({\pmb x}, {\pmb y})\big)_{{\pmb x}, {\pmb y}\in\mathcal{Z}^N}$ given by $(\ref{metro-dyn-mat})$. Define for convenience the following local empirical measures
 \begin{align*}
 \mu_j^{\iota,N}(t)=\frac{1}{N_j^{\iota}}\sum_{i\in C_j^{\iota}}\delta_{x_i(t)}, \mbox{ for $1\leq j\leq r$ and $\iota\in\{c,p\}$},
 \end{align*}
 and let $\mu^{N}(t)=\big(\mu_1^{c,N}(t), \mu_1^{p,N}(t),\ldots, \mu_r^{c,N}(t), \mu_r^{p,N}(t)\big)$ be the corresponding empirical vector. Observe that $\mu_j^{\iota,N}(t)(\cdot)=\varrho^{j,\iota}_{\cdot}({\pmb x}^N(t))$. Therefore, at any time $t\geq 0$, a given central (resp. peripheral) particle within the block $j$ jumps from $z$ to $z'$, for $z,z'\in\mathcal{Z}$, at rate $\lambda^{j,c,N}_{z,z'}\big(\mu^{N}(t),N^c_{j},N^p_{j}\big)$ given in $(\ref{rate-funct-cent-N})$ (resp. $\lambda^{j,p,N}_{z,z'}\big(\mu^{N}(t),N^c_{j},N_1^p,\ldots, N_r^p\big)$ given in $(\ref{rate-funct-periph-N})$). 
 
 Recall the important notion of multi-exchangeability for multi-class systems in the definition given below.  

 \begin{defi}
 \label{defi-multi-exch}
A sequence of random variables $(X_{n,k}, 1 \leq n \leq N_k, 1 \leq k \leq K)$ indexed by $N=(N_k)\in \mathbb{N}^K$ is said to be multi-exchangeable if its law is invariant under permutation of the indexes within the classes, that is, for $1 \leq k \leq K$ and any permutations $\sigma_k$ of $\{1,\ldots, N_k\}$, the following equality in distribution holds 
\begin{align*}
(X_{\sigma_k(n),k} , 1 \leq n \leq N_k, 1 \leq k \leq K) \overset{dist}{=} (X_{n,k}, 1 \leq n \leq N_k, 1 \leq k \leq K).
\end{align*}
 \end{defi}
 
Given the symmetry of the rate functions within each of the $2r$ classes $C_1^c,C_1^p,\ldots,C_r^c,C_r^p$, a multi-exchangeability assumption made on the initial condition ${\pmb x}^N(0)$ guarantees that, at any $t>0$, ${\pmb x}^N(t)$ is also multi-exchangeable. It follows that the local empirical measures $\mu_j^{\iota,N}(t)$ are Markov processes taking values respectively in the state spaces $\mathcal{M}^{N_j^{\iota}}_1(\mathcal{Z})$. For a detailed proof of this claim, one can consult \cite[Prop. 2.3.3]{Daw93}. The random evolution of the empirical vector $\mu^N(t)$ is summarized as follows: 

First, notice that the current setting allows almost surely at most one particle to jump at any given time. Therefore the jumps of $\mu^N(t)$ happen within each component $\mu_j^{\iota,N}(t)$ at a time and have the form $\frac{1}{N_j^{\iota}}(e_{z'}-e_z)$ for $z,z'\in\mathcal{Z}$, with $e_z$ representing the unit-coordinate vector in $\mathbb{R}^K$ in the $z$-direction. Moreover, if $\mu^N(t)=q=(q^{1,c},q^{1,p},\ldots,q^{r,c},q^{r,p})\in\prod_{j=1}^r\prod_{\iota\in\{c,p\}}\mathcal{M}^{N_j^{\iota}}_1(\mathcal{Z})$ at time $t\geq 0$, then there are $N_j^{\iota}q^{j,\iota}_z$ particles of the class $C_j^{\iota}$ at each state $z\in\mathcal{Z}$. Each of these particles jumps independently to $z'\in\mathcal{Z}$ at rate $\lambda^{j,c,N}_{z,z'}\big(q,N^{c}_{j},N^p_{j}\big)$ for $\iota=c$ and $\lambda^{j,p,N}_{z,z'}\big(q,N^c_{j},N_1^p,\ldots, N_r^p\big)$ for $\iota=p$. Therefore, the local central empirical processes $\mu_j^{c,N}(t)$ transit from $q^{j,c}$ to $q^{j,c}+\frac{1}{N_j^c}(e_{z'}-e_z)$ at rate $N_j^cq^{j,c}_z\lambda^{j,c,N}_{z,z'}\big(q,N^{c}_{j},N^p_{j}\big)$. Similarly the peripheral local empirical processes $\mu_j^{p,N}(t)$ transit from $q^{j,p}$ to $q^{j,p}+\frac{1}{N_j^p}(e_{z'}-e_z)$ at rate $N_j^pq^{j,p}_z\lambda^{j,p,N}_{z,z'}\big(q,N^c_{j},N_1^p,\ldots, N_r^p\big)$. Thence, the infinitesimal generator associated with the empirical vector $\mu^N(t)$ is given for any real-valued function $\Phi$ on $(\mathcal{M}_1(\mathcal{Z}))^{2r}$ by
\begin{align*}
\mathcal{L}^{N}\Phi(q)&=\sum_{j=1}^r\bigg[\sum_{\substack{z,z'\in\mathcal{Z}}}N_j^{c}q^{j,c}_z\lambda^{j,c,N}_{z,z'}\big(q,N_j^c,N_j^p\big)\bigg\{\Phi\bigg(q^{1,c},q^{1,p},\ldots,q^{j,c}+\frac{1}{N_j^c}(e_{z'}-e_z),\ldots,q^{r,c},q^{r,p}\bigg)-\Phi(q)\bigg\}\\
&\quad+\sum_{\substack{z,z'\in\mathcal{Z}}}N_j^{p}q^{j,p}_z\lambda^{j,p,N}_{z,z'}\big(q,N_j^c,N_1^p,\ldots,N_r^p\big)\bigg\{\Phi\bigg(q^{1,c},q^{1,p},\ldots,q^{j,p}+\frac{1}{N_j^p}(e_{z'}-e_z),\ldots,q^{r,c},q^{r,p}\bigg)-\Phi(q)\bigg\}\bigg].
\end{align*}
One might also describe the changes in the system locally by introducing the infinitesimal generators corresponding to the classes $C_j^{\iota}$ and defined for functions $f$ on $\mathcal{M}_1(\mathcal{Z})$ by
\begin{align*}
\mathcal{L}_q^{j,c,N}f(q^{j,c})&=\sum_{\substack{z,z'\in\mathcal{Z}}}N_j^{c}q^{j,c}_z\lambda^{j,c,N}_{z,z'}\big(q,N_j^c,N_j^p\big)\big[f\big(q^{j,c}+\frac{1}{N_j^c}(e_{z'}-e_z)\big)-f(q^{j,c})\big],\\
\mathcal{L}_q^{j,p,N}f(q^{j,p})&=\sum_{\substack{z,z'\in\mathcal{Z}}}N_j^{p}q^{j,p}_z\lambda^{j,p,N}_{z,z'}\big(q,N_j^c,N_1^p,\ldots,N_r^p\big)\big[f\big(q^{j,p}+\frac{1}{N_j^p}(e_{z'}-e_z)\big)-f(q^{j,p})\big].
\end{align*}
 
As mentioned in the introduction, one classical question in the study of interacting particle systems is their large-scale behavior, namely, their behavior when the total number of particles goes to infinity. In the next result, we state that a law of large numbers holds for the empirical vector $\mu^N(t)$ as $N\rightarrow\infty$. It is worthwhile to mention that contrary to the models studied in \cite{Daw+Sid+zha2020, Daw+Sid+zha2021}, the rate functions here are convergent.  
 
\begin{prop} 
  Suppose that the initial condition $\mu^{N}(0)=\big(\mu_1^{c,N}(0), \mu_1^{p,N}(0),\ldots, \mu_r^{c,N}(0), \mu_r^{p,N}(0)\big)$ converges weakly, as $N\rightarrow\infty$, towards $\nu=(\nu_{1}^{c,N},\nu_{1}^{p,N},\ldots,\nu_{r}^{c,N},\nu_{r}^{p,N})\in(\mathcal{M}_1(\mathcal{Z}))^{2r}$. Then, the empirical vector process $\mu^N(t)$ converges in probability and uniformly over  compact time intervals toward the solution $\mu$ to the following McKean-Vlasov system 
\begin{equation}
\begin{split}
\left\{ \begin{array}{lcl}\dot{\mu}_j^c(t)=\mu_j^c(t)A^{j,c}_{\mu (t)}, & &  \\
 \dot{\mu}_j^p(t)=\mu_j^p(t)A^{j,p}_{\mu(t)},  & & \\
 \mu_j^c(0)=\nu_j^c,\mu_j^p(0)=\nu_j^p, & & \\
 1\leq j\leq r.  & &
 \end{array}\right.
\end{split}
\label{McKean-Vlas-syst}
\end{equation}
\label{LLN-prop}
\end{prop}

\proof Since $\mu^N(t)$ is a pure jump Markov process, one might rely on the classical Kurtz theorem \cite{Kurtz70}. We first verify the conditions of its application. Denote $E=\mathcal{M}_1(\mathcal{Z})$ and $E_{N_j^{\iota}}=\mathcal{M}^{N_j^{\iota}}_1(\mathcal{Z})$ for all $1\leq j\leq r$ and $\iota\in\{c,p\}$. Define the following functions 
\begin{align*}
F_q^{j,c,N}(q^{j,c})&=\mathcal{L}_q^{j,c,N}q^{j,c}=\sum_{\substack{z,z'\in\mathcal{Z}}}N_j^{c}q^{j,c}_z\lambda^{j,c,N}_{z,z'}\big(q,N_j^c,N_j^p\big)\frac{1}{N_j^c}(e_{z'}-e_z),\\
F_q^{j,p,N}(q^{j,p})&=\mathcal{L}_q^{j,p,N}q^{j,p}=\sum_{\substack{z,z'\in\mathcal{Z}}}N_j^{p}q^{j,p}_z\lambda^{j,p,N}_{z,z'}\big(q,N_j^c,N_1^p,\ldots,N_r^p\big)\frac{1}{N_j^p}(e_{z'}-e_z),
\end{align*}
\begin{align*}
F_q^{j,c}(q^{j,c})&=\sum_{\substack{z,z'\in\mathcal{Z}}}q^{j,c}_z\lambda^{j,c}_{z,z'}\big(q\big)(e_{z'}-e_z),\\
F_q^{j,p}(q^{j,p})&=\sum_{\substack{z,z'\in\mathcal{Z}}}q^{j,p}_z\lambda^{j,p}_{z,z'}\big(q\big)(e_{z'}-e_z),
\end{align*}
for $1\leq j\leq r$ and $q=(q^{1,c},q^{1,p},\ldots,q^{r,c},q^{r,p})\in\prod_{j=1}^r\prod_{\iota\in\{c,p\}} E_{N_j^{\iota}}$. Observe from $(\ref{lambd-gibbs})$ that the rate functions $\lambda^{j,\iota}_{z,z'}$ are Lipschitz in $q^{j,\iota}$ for all $1\leq j\leq r$ and $\iota\in\{c,p\}$. Therefore, since $q^{j,\iota}$ are probability measures, one can find a constant $M$ such that, for any $\bar{q}=(q^{1,c},q^{1,p},\ldots,\bar{q}^{r,\iota},\ldots,q^{r,c},q^{r,p})\in\prod_{j=1}^r\prod_{\iota\in\{c,p\}}E_{N_j^{\iota}}$,
\begin{align}
\|F_{q}^{j,\iota}(q^{j,\iota})-F_{\bar{q}}^{j,\iota}(\bar{q}^{j,\iota})\|\leq M \|q^{j,\iota}-\bar{q}^{j,\iota}\|,
\label{Lipschi}
\end{align}
where $\|\cdot\|$ denotes a distance on $\mathbb{R}^{K}$. In addition, recall that the functions $\lambda^{j,\iota,N}_{z,z'}$ convergence toward $\lambda^{j,\iota}_{z,z'}$ as $N\rightarrow\infty$. Hence, it is easy to see that
\begin{align*}
\lim_{N\rightarrow\infty}\sup_{q^{j,\iota}\in E_{N_j^{\iota}} \cap E} |F_q^{j,\iota,N}(q^{j,\iota})-F_q^{j,\iota}(q^{j,\iota})|=0, \mbox{  for all $1\leq j\leq r$ and $\iota\in\{c,p\}$. } 
\end{align*}

Now, straightforward calculations allow us to verify that the $z$-th component of the $K$-dimensional vector $F_q^{j,\iota}(q^{j,\iota})$ is equal to $\sum_{z'\neq z}q^{j,\iota}_{z'}\lambda^{j,\iota}_{z',z}\big(q\big)-\sum_{z'\neq z}q^{j,\iota}_{z}\lambda^{j,\iota}_{z,z'}\big(q\big)=\sum_{z'\in\mathcal{Z}}q^{j,\iota}_{z'}\lambda^{j,\iota}_{z',z}\big(q\big)$, which in turn corresponds to the $z$-th component of the raw vector $q^{j,\iota}A^{j,\iota}_{q}$. Finally, notice that 
\begin{align*}
\sup_{N_j^c}\sup_{q^{j,c}\in E_{N_j^c}}\sum_{\substack{z,z'\in\mathcal{Z}}}N_j^{c}q^{j,c}_z\lambda^{j,c,N}_{z,z'}\big(q,N_j^c,N_j^p\big)\frac{1}{N_j^c}\|e_{z'}-e_z\|
\end{align*}
and 
\begin{align*}
\sup_{N_j^p}\sup_{q^{j,p}\in E_{N_j^p}}\sum_{\substack{z,z'\in\mathcal{Z}}}N_j^{p}q^{j,p}_z\lambda^{j,p,N}_{z,z'}\big(q,N_j^c,N_1^p,\ldots,N_r^p\big)\frac{1}{N_j^p}\|e_{z'}-e_z\|
\end{align*}
are bounded given the form of the rate functions $\lambda^{j,c,N}_{z,z'}$ and $\lambda^{j,p,N}_{z,z'}$ in $(\ref{rate-funct-cent-N})$ and $(\ref{rate-funct-periph-N})$ and the fact that $\mathcal{Z}$ is finite. Thus condition $(2.9)$ in \cite{Kurtz70} is verified. To verify condition $(2.10)$ of \cite{Kurtz70}, define $\varepsilon_{N_j^{\iota}}=\frac{C}{N_j^{\iota}}$ where $C>\|e_{z'}-e_z\|$. Thus $\varepsilon_{N_j^c}\downarrow 0$ as $N\rightarrow\infty$, and the following converges hold true
\begin{align*}
\lim_{N_j^c\rightarrow\infty}\sup_{q^{j,c}\in E_{N_j^c}}\sum_{\substack{z,z'\\\frac{1}{N_j^c}\|e_{z'}-e_z\|>\varepsilon_{N_j^c}}}N_j^{c}q^{j,c}_z\lambda^{j,c,N}_{z,z'}\big(q,N_j^c,N_j^p\big)\frac{1}{N_j^c}\|e_{z'}-e_z\|=0,
\end{align*}
\begin{align*}
\lim_{N_j^p\rightarrow\infty}\sup_{q^{j,p}\in E_{N_j^p}}\sum_{\substack{z,z'\\\frac{1}{N_j^p}\|e_{z'}-e_z\|>\varepsilon_{N_j^p}}}N_j^{p}q^{j,p}_z\lambda^{j,p,N}_{z,z'}\big(q,N_j^c,N_1^p,\ldots,N_r^p\big)\frac{1}{N_j^p}\|e_{z'}-e_z\|=0.
\end{align*}
Thus condition $(2.10)$ in \cite{Kurtz70} holds. Therefore, we are now in the position to apply \cite[Theorem. 2.11]{Kurtz70} since all the related conditions are satisfied. Thus each $\mu_j^{\iota,N}(t)$ converges in probability and uniformly over any time interval $[0,T]$ towards the solution $\mu_j^{\iota}(t)$ of the differential equation $\frac{d}{dt}\mu_j^{\iota}(t)=\mu_j^{\iota}(t)A^{j,\iota}_{\mu(t)}$. Define 
\begin{align*}
F(q)=\bigg(F_q^{1,c}(q^{1,c}),F_q^{1,p}(q^{1,p}),\ldots,F_q^{r,c}(q^{r,c}),F_q^{r,p}(q^{r,p})   \bigg).
\end{align*}
From $(\ref{Lipschi})$, it follows that $F$ is Lipschitz. Therefore standard arguments show that the differential equation $\dot{\mu}(t)=F(\mu(t))$ has a unique solution. This concludes the proof.
\carre

\section{Stability of the McKean-Vlasov system}
\label{stab-sect}
The law of large numbers established in the last section characterizes the large-scale behavior of the $N$-particles system over finite time intervals. In particular, as $N\rightarrow\infty$, the empirical vector $\mu^N(t)$ converges in probability towards the solution $\mu(t)$ to the McKean-Vlasov system $(\ref{McKean-Vlas-syst})$. Thence, when the total number of particles in the system is very large, one might approximate the behavior of $\mu^N(t)$ by $\mu(t)$ over finite time intervals. However, when time $t\rightarrow\infty$, this approximation might no longer be accurate. In particular, this will depend on the uniqueness of equilibrium points of the McKean-Vlasov system and their stability. In the case where there are multiple equilibria or a unique unstable equilibrium, metastable phenomena can arise and care must be taken in using the approximation. One can consult \cite{Bork+Sund2012, Daw+Sid+zha2021} for a detailed discussion. The long-time behavior of the large $N$-particles system is thus related to the stability of the McKean-Vlasov system $(\ref{McKean-Vlas-syst})$. The goal of the current section is to investigate it by constructing a suitable Lyapunov function. The approach we take is based on the calculation of the limit of suitably normalized relative entropies. This idea was introduced in \cite{Budh+Du(a)2015, Budh+Du(b)2015} to study the stability of Kolmogorov forward equations arising as the limit of mean-field systems with jumps on complete graphs. We thus generalize this method to multi-class mean-field systems with jumps. 

Let us first introduce some definitions. Recall that $|\mathcal{Z}|=K$ is the number of possible states. Let $\mathcal{S}=\{ p\in\mathbb{R}^K: p_z\geq 0, z \in\mathcal{Z}, \sum_{z\in\mathcal{Z}} p_z = 1 \}$ be the $(K-1)$-dimensional simplex, and let $\mathcal{S}^{\circ}=\{ p\in \mathcal{S}: p_z> 0,\mbox{ for all $z \in\mathcal{Z}$}\}$ denotes its relative interior. Notice that the space $\mathcal{M}_1(\mathcal{Z})$ of probability measures on $\mathcal{Z}$ can be identified with $\mathcal{S}$. Thus, the empirical vector $(\varrho(t),t\geq 0)$ takes values in $\mathcal{S}^{2r}$.

\begin{defi}
A point $\phi=(\phi_j^c,\phi_j^p,1\leq j\leq r)\in\mathcal{S}^{2r} $  is said to be a fixed point of the McKean-Vlasov system $(\ref{McKean-Vlas-syst})$ if the right-hand side of $(\ref{McKean-Vlas-syst})$ evaluated at $\mu=\phi$ is equal to zero, namely,
\begin{equation}
\begin{split}
\left\{ \begin{array}{lcl}\phi_j^cA_{\phi}=0, & &  \\
                          \phi_j^pA_{\phi}=0,  & & \\
 1\leq j\leq r. & &
 \end{array}\right.
\end{split}
\label{stat-point}
\end{equation}
\end{defi} 

\begin{defi}
A fixed point $\phi=(\phi_j^c,\phi_j^p,1\leq j\leq r)\in \mathcal{S^{\circ}}^{2r}$ of the McKean-Vlasov system $(\ref{McKean-Vlas-syst})$ is said to be locally stable if there exists a relatively open subset $\Gamma$ of $\mathcal{S}^{2r}$ that contains $\phi$ and has the property that whenever $\mu(0)\in\Gamma$, the solution $\mu(t)$ to $(\ref{McKean-Vlas-syst})$ with initial condition $\mu(0)$ converges to $\phi$ as $t\rightarrow\infty$.
\end{defi}
 
\begin{defi}
Let $\phi=(\phi_j^c,\phi_j^p,1\leq j\leq r)\in\mathcal{S^{\circ}}^{2r}$  be a fixed point of $(\ref{McKean-Vlas-syst})$ and let $\Gamma$ be a relatively open subset of $\mathcal{S}^{2r}$ that contains $\phi$. A function $J: \Gamma\rightarrow\mathbb{R}$ is called positive definite if for some $K^*\in\mathbb{R}$, the sets $M_K= \{r \in\bar{\Gamma}: J(r) \leq K\}$ decrease continuously to ${\phi}$ as $K\downarrow K^*$.
\end{defi}
  
\begin{defi}
Let $\phi\in \mathcal{S^{\circ}}^{2r}$ be a fixed point of the McKean-Vlasov system $(\ref{McKean-Vlas-syst})$, and let $\Gamma$ be a relatively open subset of $\mathcal{S}^{2r}$ that contains $\phi$. A positive definite, $\mathcal{C}_1$ and uniformly continuous function $J: \Gamma\rightarrow\mathbb{R}$ is said to be a local Lyapunov function associated with $(\Gamma,\phi)$ for the McKean-Vlasov system  $(\ref{McKean-Vlas-syst})$ if, given any $\mu(0)\in\Gamma$, the solution $\mu(t)$ to $(\ref{McKean-Vlas-syst})$  with initial condition $\mu(0)$ satisfies $\frac{d}{dt}J(\mu(t))<0$ for all $0\leq t < \tau$ such that $\mu(t)\neq \phi$, where $\tau=\inf\{t \geq 0 : \mu(t) \in \Gamma^c\}$. Moreover, if $\Gamma= {\mathcal{S}^{\circ}}^{2r} $, we say that $J$ is a Lyapunov function.
\label{lyapu-def}
\end{defi} 

 The next classical result shows that the existence of a local Lyapunov function is equivalent to the local stability of an equilibrium of the McKean-Vlasov system $(\ref{McKean-Vlas-syst})$. We give detailed proof tailored for the specific system $(\ref{McKean-Vlas-syst})$. 
 
 \begin{prop}
 Let $\phi\in \mathcal{S^{\circ}}^{2r}$ be a fixed point of the McKean-Vlasov system $(\ref{McKean-Vlas-syst})$, and let $\Gamma$ be some relatively open subset of $\mathcal{S}^{2r}$ that contains $\phi$. Suppose that there exists a local Lyapunov function associated with $(\Gamma, \phi)$. Then $\phi$ is locally stable.
 \label{stab-lyap}
 \end{prop} 

\proof Let $J$ be a Lyapunov function associated with $(\Gamma, \phi)$. Therefore, $J$ is positive definite. Thus, for some $K^*\in\mathbb{R}$, the sets $M_K= \{r \in\bar{\Gamma}: J(r) \leq K\}$ decrease continuously to ${\phi}$ as $K\downarrow K^*$. One can then find some $L>K^*$ and an open set $\mathcal{O}$ such that $\phi\in\mathcal{O}\subset M_L\subset\Gamma$. Let $\mu(0)\in \mathcal{O}$ such that $\mu(0)\neq\phi$. By the chain rule property together with $(\ref{McKean-Vlas-syst})$ we find
\begin{align}
\frac{d}{dt}J(\mu(t))=\bigg\langle D J \big(\mu(t)\big), \big(\mu_1^c(t)A^{1,c}_{\mu (t)}, \mu_1^p(t)A^{1,p}_{\mu (t)},\ldots,\mu_r^c(t)A^{r,c}_{\mu (t)},\mu_r^p(t)A^{r,p}_{\mu (t)}\big) \bigg\rangle,
\label{deriv-lyap}
\end{align}
where $DJ(\mu(t))$ is the derivative (gradient's transpose) of $J$ at the solution $\mu(t)$ of the McKean-Vlasov system $(\ref{McKean-Vlas-syst})$, and $\langle\cdot,\cdot\rangle$ denotes the scalar product on $\mathbb{R}^{2r}$. Define the stopping time $\tau=\inf\{t>0:\mu(t)\in \Gamma^c\}$. Suppose that $\tau<\infty$. By definition of Lyapunov function we have that $\frac{d}{dt}J(\mu(t))<0$ for all $0\leq t \leq \tau$ and $\mu(t)\neq\phi$. Moreover, observe that the function $q\rightarrow \big\langle D J \big(q\big), \big(q_1^cA^{1,c}_{q}, q_1^pA^{1,p}_{q},\ldots,q_r^cA^{r,c}_{\mu },q_r^pA^{r,p}_{q}\big) \big\rangle$ is continuous on $\Gamma$. Then we have necessarily $\frac{d}{dt}J(\mu(t))=0$ for $\mu(t)=\phi$ since $\frac{d}{dt}J(\mu(t))<0$ elsewhere on $\Gamma$ (no discontinuity). Thus $\frac{d}{dt}J(\mu(t))\leq 0$ for all $0\leq t\leq\tau$. Furthermore, since $\mathcal{S}^{2r}$ is a metric space, $\mathbb{R}$ is a complete metric space, and $J$ is uniformly continuous, we have that the function $J$ can be extended continuously on the closure $\bar{\Gamma}$ of $\Gamma$. Using this together with $\frac{d}{dt}J(\mu(t))\leq 0$ for all $0\leq t \leq \tau$ and since $\mu(0)\in M_L$ proves that $J(\mu(\tau)) \leq L$ (since $J$ is decreasing when $t\in [0,\tau]$). We thus deduce that $\mu(\tau)\subset M_L\subset\Gamma$. This means that $\tau=\infty$ and hence again by definition of Lyapunov function
\begin{align}
\mbox{   $\frac{d}{dt}J(\mu(t))< 0$ for all $t\geq 0$ with $\mu(t)\neq\phi$.}
\label{neg-deriv-lyap}
\end{align}

   Let $ K^*< K_n<L$ for $n\geq 1$  be a decreasing sequence of real numbers going to $K^*$, and define the corresponding stopping times  
   \begin{align*}
   \tau_n=\inf\{t>0:\mu(t)\in M_{K_n}\}. 
   \end{align*}
We next prove that $\tau_n<\infty$ for all $n\geq 1$. First, if $\mu(0)\in M_{K_n}$, this follows immediately by the above arguments. Suppose that $\mu(0)\notin M_{K_1}$. Then, by the positive definitness of the function $J$, we can always find $\varepsilon>0$ such that the ball $B(\phi,\varepsilon)$ centered at $\phi$ with radius $\varepsilon$ satisfies $B(\phi,\varepsilon)\cap\mathcal{S}^{2r}\subset M_{K_1}\subset M_L$. Let $q\in (B(\phi,\varepsilon))^c\cap M_L$ and define by $\mu^q(t)$ the solution to the McKean-Vlasov system $(\ref{McKean-Vlas-syst})$ with initial condition $\mu^q(0)=q$. Therefore, using $(\ref{neg-deriv-lyap})$, one obtains that $\frac{d}{dt}J(\mu^q(t))\big|_{t=0}< 0$.  Recall that the function $q\rightarrow \big\langle D J \big(q\big), \big(q_1^cA^{1,c}_{q}, q_1^pA^{1,p}_{q},\ldots,q_r^cA^{r,c}_{\mu },q_r^pA^{r,p}_{q}\big) \big\rangle$ is continuous on $\Gamma$ , and that $B(\phi,\varepsilon)\cap M_L$ is a closed subset of $\Gamma$, thus   
   \begin{align*}
   \sup_{q\in ( B(\phi,\varepsilon))^c\cap M_L}\big\langle D J \big(q\big), \big(q_1^cA^{1,c}_{q}, q_1^pA^{1,p}_{q},\ldots,q_r^cA^{r,c}_{\mu },q_r^pA^{r,p}_{q}\big) \big\rangle<0.
\end{align*}    
In addition, since $ B(\phi,\varepsilon)\cap\mathcal{S}^{2r}\subset M_{K_1}$, we have that $\mu(t)\in (B(\phi,\varepsilon))^c\cap M_L$ for all $t\leq\tau_1$ from which we deduce that $\sup_{t\leq\tau_1}\frac{d}{dt}J(\mu(t))\leq 0$. Then $\tau_1<\infty$. Repeating the same argument for all $n>1$ gives that $\tau_n<\infty$ for all $n\geq 1$. This concludes the proof. \carre

\subsection{Limit of relative entropies}

We now construct a candidate Lyapunov function for the McKean-Vlasov system $(\ref{McKean-Vlas-syst})$ as the limit of suitably scaled relative entropies. We start by recalling the form of the relative entropy function denoted by $R(\cdot \|\cdot)$ and defined, for all $p\in(\mathcal{M}_1(\mathcal{Z}))^{2r}$, by
\begin{align}
R(p\|q)=\sum_{z\in\mathcal{Z}^N}p_z\log\left(\frac{p_z}{q_z}\right).
\label{relat-entro-def}
\end{align} 
Notice that the relative entropy function plays an important role in various fields of mathematics. For an account of its properties and applications, one can consult e.g. \cite{Cover+Thomas2006}. Define the function $\bar{F}^N$ given, for any $q=(q^{1,c},q^{1,p},\ldots,q^{r,c},q^{r,p})\in\big(\mathcal{M}_1(\mathcal{Z})\big)^{2r}$, by
\begin{align}
\bar{F}^N(q)=\frac{1}{N}R\bigg(\otimes^{N_1^c}q^{1,c}\otimes^{N_1^p}q^{1,p}\cdots\otimes^{N_r^c}q^{r,c}\otimes^{N_r^p}q^{r,p}\big\|\pi^N({\pmb x})\bigg),
\label{Lyap-func}
\end{align}
where $\pi^N({\pmb x})$ is the Gibbs measure introduced in $(\ref{stat-distr-gibbs})$. The idea now is to identify $\lim_{N\rightarrow\infty}\bar{F}^N(q)$ and then investigate its Lyapunov function properties for the McKean-Vlasov limiting system $(\ref{McKean-Vlas-syst})$. The calculation of this limit relies on the following Laplace principle for empirical vectors which is an extension of the Sanov's theorem for empirical measures given in \cite[Th. $2.2.1$]{Dupuis+Ellis97}. 

\begin{prop}
\label{Sanov-vect-meas-prop}
Let $\mathcal{Z}$ be a Polish space and let $\nu$ be a probability measure on $\mathcal{Z}$. Let $r\geq 1$, and let $\{x^{j}_i\}_{1\leq i\leq N_j}$, for $1\leq j\leq r$, be $r$ sequences of $\mathcal{Z}$-valued i.i.d random variables with a common distribution $\nu$. Suppose that $\sum_{j=1}^rN_j=N$ and that there exist positive reals $\alpha_1,\ldots,\alpha_r$ such that $\lim_{N\rightarrow\infty}\frac{N_j}{N}=\alpha_j$ and $\sum_j\alpha_j=1$. Denote by $\mu^N=\big(\mu^N_{1},\ldots,\mu^N_{r}\big)$ the empirical vector associated with these sequences, i.e., for each $1\leq j\leq r$, $\mu^N_{j}=\frac{1}{N_j}\sum_{i=1}^{N_j}\delta_{x^{j}_i}$. Then, for all sequence $\{h^N,N\in\mathbb{N}\}$ of bounded continuous functions mapping $(\mathcal{M}_1(\mathcal{Z}))^r$ into $\mathbb{R}$ and  converging to some $h$ as $N\rightarrow\infty$, the following Laplace principle holds
\begin{align*}
\lim_{N\rightarrow\infty}\frac{1}{N}\log\mathbb{E}\bigg[\exp\bigg\{-Nh^N(\mu^N)\bigg\}\bigg]=-\inf_{\gamma\in(\mathcal{M}_1(\mathcal{Z}))^{r}} \bigg\{R\bigg(\sum_{j=1}^r\alpha_j\gamma_j\big\|\nu\bigg)+h(\gamma)\bigg\}.
\end{align*}  
\end{prop}

\proof We follow the strategy elaborated in \cite{Dupuis+Ellis97}. We first establish a variational formulation, then we prove the convergence by establishing lower and upper bounds.  
 
\paragraph{Step 1: Variational formulation.} Set
\begin{align*}
W^N&=-\frac{1}{N}\log\mathbb{E}\bigg[\exp\bigg\{-Nh^N(\mu^N)\bigg\}\bigg]\\
   &=-\frac{1}{N}\log\int_{\prod_{j=1}^r\mathcal{Z}^{N_j}} \exp\{-Nh^N(\mu^N)\}d\big(\otimes_{j=1}^r\otimes^{N_j}\nu\big).
\end{align*}
Therefore, one can represent $W^N$ by the following variational formula \cite[Prop. 1.4.2]{Dupuis+Ellis97}
\begin{align*}
W^N=\inf_{\eta^N\in\mathcal{M}_1\big(\prod_{j=1}^r\mathcal{Z}^{N_j}\big)}\bigg\{\frac{1}{N} R\bigg(\eta^N\big\|\otimes_{j=1}^r\otimes^{N_j}\nu\bigg)+ \int_{\prod_{j=1}^r\mathcal{Z}^{N_j}}h^N(\mu^N) d\eta^N\bigg\}.
\end{align*} 
Using the decomposition property of probability measures on product spaces, one can observe that  
\begin{align*}
\eta^N\big((dx_1^{1}\times\ldots\times x_{N_1}^{1})\times\ldots \times (dx_1^{r}\times\ldots\times dx_{N_r}^{r})\big)
\end{align*}  
is equivalent to   
\begin{align*}
\eta_1^N\big(dx_1^{1}\times\ldots\times x_{N_1}^{1}\big)\otimes\eta_2^N\big((dx_1^{2}\times&\ldots\times dx_{N_2}^{2})|(x_1^{1},\ldots, x_{N_1}^{1})\big)\otimes\cdots\\
&\cdots\otimes \eta^N_r \big((dx_1^{r}\times\ldots\times dx_{N_r}^{r})| (x_1^{1},\ldots, x_{N_1}^{1}),\ldots,(x_1^{r-1},\ldots, x_{N_{r-1}}^{r-1})\big),
\end{align*}
where $\eta_1^N$ is the projection onto the space $\mathcal{Z}^{N_1}$ and, for $j\geq 2$, $\eta^N_j$ is the conditional distribution on $\mathcal{Z}^{N_j}$ given $(x_1^{1},\ldots, x_{N_1}^{1}),\ldots,(x_1^{j-1},\ldots, x_{N_1}^{j-1})$. Similarly, the last display can be further factorized as  
\begin{align*}
\eta_{11}^N\big(dx_1^{1}\big)\otimes\ldots&\otimes\eta^N_{1N^1}\big(dx_{N_1}^{1}|x_{1}^{1},\ldots,x_{N_1-1}^{1}\big)\otimes\eta_{21}^N\big(dx_1^{2}|x_1^{1},\ldots, x_{N_1}^{1}\big)\otimes\eta_{22}^N\big(dx_2^{2}|x_1^{1},\ldots, x_{N_1}^{1},x_1^{2}\big)\otimes\cdots\\ 
&\cdots\otimes \eta_{2N_2}^N\big(dx_{N_2}^{2}|x_1^{1},\ldots, x_{N_1}^{1},x_1^{2},\ldots, x_{N_2-1}^{2}\big)\otimes
\cdots\otimes \eta^N_{rN_r} \big(dx_{N_r}^{r}| x_1^{1},\ldots, x_{N_1}^{1},\ldots,x_1^{r},\ldots, x_{N_r-1}^{r}\big),
\end{align*}
where, for $1\leq j\leq r$  and $ 1\leq i \leq N_j$, $\eta_{ji}$ is the conditional distribution on $\mathcal{Z}$ given $x_1^{1},\ldots, x_{N_1}^{1},\ldots,x_1^{j},\ldots, x_{i-1}^{j}$. Using this together with the chain rule property of relative entropy (see \cite[Th. 2.5.3]{Cover+Thomas2006} or \cite[Th. C.3.1]{Dupuis+Ellis97}) one obtains   
 \begin{align*}
 R\big(\eta^N\|\nu^N\big)=\int_{\prod_{j=1}^r\mathcal{Z}^{N_j}} \sum_{j=1}^r\sum_{i=1}^{N_j} R\bigg(\eta_{ji}^N\big(\cdot |x_1^{1},\ldots, x_{N_1}^{1},\ldots,x_1^{j},\ldots, x_{i-1}^{j}\big)\big\| \nu (\cdot)\bigg) d\eta^N, 
 \end{align*}
from which we deduce the new variational formulation 
 \begin{align*}
W^N=\inf_{\eta^N\in\mathcal{M}_1(\prod\limits_{j=1}^r\mathcal{Z}^{N_j})}\bigg\{\frac{1}{N} \int_{\prod\limits_{j=1}^r\mathcal{Z}^{N_j}} \sum_{j=1}^r\sum_{i=1}^{N_j}& R\bigg(\eta_{ji}^N(\cdot |x_1^{1},\ldots, x_{N_1}^{1},\ldots,x_1^{j},\ldots, x_{i-1}^{j})\big\| \nu (\cdot)\bigg)  d\eta^N \\
&+ \int_{\prod\limits_{j=1}^r\mathcal{Z}^{N_j}}h^N(\mu^N) d\eta^N\bigg\}.\\
\end{align*} 
 
Notice that the infimum in the last formulation can be replaced by the infimum  over all sequences $\{\eta^N_{ji}\}$ of conditional distributions  on $\mathcal{Z}$ given $x_1^{1},\ldots, x_{N_1}^{1},\ldots,x_1^{j},\ldots, x_{i-1}^{j}$. Now, given a sequence of conditional distributions $\{\eta^N_{ji}\}$, define the sequence of $\mathcal{Z}$-valued random variables $\{\bar{x}^{j}_i\}_{1\leq i\leq N_j}$ by specifying their distributions recursively as follows: 
\begin{align*}
\mathbb{P}(\bar{x}^{j}_i\in dy|\bar{x}_1^{1},\ldots, \bar{x}_{N_1}^{1},\ldots,\bar{x}_1^{j},\ldots, \bar{x}_{i-1}^{j})=\eta_{ji}^N\big(dy |\bar{x}_1^{1},\ldots, \bar{x}_{N_1}^{1},\ldots,\bar{x}_1^{j},\ldots, \bar{x}_{i-1}^{j}\big),
\end{align*}
and set the corresponding empirical measure vector to $\bar{\mu}^N=\big(\bar{\mu}^N_{1},\ldots,\bar{\mu}^N_{r}\big)$ where, for $1\leq j\leq r$, 
\begin{align*}
\bar{\mu}^N_{j}=\frac{1}{N_j}\sum_{i=1}^{N_j}\delta_{\bar{x}^{j}_i}.
\end{align*}
Thence, one can easily verify that the following variational formulation holds true
 \begin{align}
W^N&=\inf_{\{\eta_{ji}^N\}}\bar{\mathbb{E}}\bigg\{\frac{1}{N} \sum_{j=1}^r\sum_{i=1}^{N_j} R\big(\eta_{ji}^N\big(\cdot |\bar{x}_1^{1},\ldots, \bar{x}_{N_1}^{1},\ldots,\bar{x}_1^{j},\ldots, \bar{x}_{i-1}^{j}\big)\big\| \nu (\cdot)\big) + h^N(\bar{\mu}^N) \bigg\},
\label{variat-form}
\end{align} 
where $ \bar{\mathbb{E}}$ denotes the expectation with respect to the random variables $\{\bar{x}^{j}_i: 1\leq i\leq N_j, 1\leq j\leq r\}$.

\paragraph{Step 2: Upper bound.} Let $\{\gamma_j\}_{1\leq j\leq r}\in \mathcal{M}_1(\mathcal{Z})$ be a sequence of probability measures on $\mathcal{Z}$ and set $\eta^N_{ji}=\gamma_j$ for all $1\leq i\leq N_j$. Then, by $(\ref{variat-form})$, one obtains the following upper bound
\begin{align*}
W^N\leq \mathbb{E}\bigg\{R\bigg(\sum_{j=1}^r\frac{N_j}{N}\gamma_j (\cdot)\big\| \nu (\cdot)\bigg)+ h^N(\bar{\mu}^N) \bigg\}.
\end{align*}
Moreover, since $\eta^N_{ji}=\gamma_j$ for all $1\leq i\leq N_j$, the random variables $\bar{x}^{j}_i$, for $1\leq i\leq N_j$, are i.i.d. Furthermore, since the measurable functions $h^N$ are bounded continuous and convergent towards $h$, using the dominated convergence theorem, the laws of large numbers, and the convergence of the proportions $\frac{N_j}{N}$ towards $\alpha_j$, one obtains 
 \begin{align*}
 \limsup_{N\rightarrow\infty} W^N\leq   R\bigg(\sum_{j=1}^r\alpha_j\gamma_j (\cdot)\big\| \nu (\cdot)\bigg)+ h(\gamma), 
 \end{align*}
 with $\gamma=(\gamma_1,\ldots,\gamma_r)$. Therefore, since the measures $\{\gamma_j, 1\leq j \leq r\}$ are arbitrary chosen, the following upper bound holds true
  \begin{align}
 \limsup_{N\rightarrow\infty} W^N\leq  \inf_{\gamma\in(\mathcal{M}_1(\mathcal{Z}))^r}\bigg\{ R\bigg(\sum_{j=1}^r\alpha_j\gamma_j(\cdot)\big\| \nu (\cdot)\bigg)+ h(\gamma) \bigg\}.
\label{upp-bound-variat-form} 
 \end{align}

 \paragraph{Step 3: Lower bound.}  First, by the convexity of relative entropy and the Jensen's inequality, the following upper bound holds true 
 \begin{align}
W^N\geq \inf_{\{\eta_{ji}^N\}}\bar{\mathbb{E}}\bigg\{R\bigg(\frac{1}{N} \sum_{j=1}^r\sum_{i=1}^{N_j} \eta_{ji}^N\big(\cdot |\bar{x}_1^{1},\ldots, \bar{x}_{N_1}^{1},\ldots,\bar{x}_1^{j},\ldots, \bar{x}_{i-1}^{j}\big)\bigg\| \nu (\cdot)\bigg) + h^N(\bar{\mu}^N) \bigg\}.
\end{align} 
Set, for all $1\leq j\leq r$,
 \begin{align*}
\chi_j^N(\cdot)= \frac{1}{N_j} \sum_{i=1}^{N_j} \eta_{ji}^N\big(\cdot |\bar{x}_1^{1},\ldots, \bar{x}_{N_1}^{1},\ldots,\bar{x}_1^{j},\ldots, \bar{x}_{i-1}^{j}\big).
 \end{align*} 
 Therefore, for any $\varepsilon>0$, one can always find a sequence $\{\eta_{ji}^N\}$ such that 
  \begin{align*}
W^N+\varepsilon\geq \bar{\mathbb{E}}\bigg[ R\bigg( \sum_{j=1}^r \frac{N_j}{N} \chi_j^N(\cdot)\big\| \nu (\cdot)\bigg)+ h^N(\bar{\mu}^N) \bigg].
\end{align*}
 
Now, for all $1\leq j\leq r$ and $1\leq i\leq N_j$, let $\mathcal{F}_{j,i}$ be the $\sigma$-algebra generated by the random variables $\bar{x}_1^{1},\ldots, \bar{x}_{N_1}^{1},\ldots,\bar{x}_1^{j},\ldots, \bar{x}_{i}^{j}$. Then, for any bounded measurable function $g:\mathcal{Z}\rightarrow\mathbb{R}$, it is easy to verify that
\begin{align*}
\bar{\mathbb{E}}\bigg[g(\bar{x}_{i}^{j})-\int_{\mathcal{Z}}g(y)\eta_{ji}^N\big(dy |\bar{x}_1^{1},\ldots, \bar{x}_{N_1}^{1},\ldots,\bar{x}_1^{j},\ldots, \bar{x}_{i-1}^{j}\big)\bigg|\mathcal{F}_{j,i-1}\bigg]=0.
\end{align*}  
One thus deduces that
\begin{align*}
\bigg\{g(\bar{x}_{i}^{j})-\int_{\mathcal{Z}}g(y)\eta_{ji}^N\big(dy |\bar{x}_1^{1},\ldots, \bar{x}_{N_1}^{1},\ldots,\bar{x}_1^{j},\ldots, \bar{x}_{i-1}^{j}\big)\bigg\}
\end{align*} 
 forms a martingale difference sequence with respect to $\mathcal{F}_{j,i}$. In addition, straightforward calculations give, for all $1\leq j\leq r$,
 \begin{align*}
 \int_{\mathcal{Z}}g(y)d\bar{\mu}^N_j- \int_{\mathcal{Z}}g(y)d\chi_j^N=\frac{1}{N_j}\sum_{i=1}^{N_j}\bigg(g(\bar{x}_i^j)-\int_{\mathcal{Z}}g(y)\eta_{ji}^N\big(dy |\bar{x}_1^{1},\ldots, \bar{x}_{N_1}^{1},\ldots,\bar{x}_1^{j},\ldots, \bar{x}_{i-1}^{j}\big)\bigg).
 \end{align*}
 Hence, recalling that the terms of a martingale difference sequence are uncorrelated we obtain by the Markov's inequality that, for any $\varepsilon>0$,
 \begin{align*}
\mathbb{P}\bigg( \bigg|\int_{\mathcal{Z}}g(y)d\bar{\mu}^N_j- \int_{\mathcal{Z}}g(y)d\chi_j^N\bigg|>\varepsilon\bigg) &\leq \frac{1}{\varepsilon^2 N^2_j}\bar{\mathbb{E}} \bigg[ \bigg| \sum_{i=1}^{N_j}\bigg(g(\bar{x}_j^j)-\int_{\mathcal{Z}}g(y)\eta_{ji}^N\big(dy |\bar{x}_1^{1},\ldots, \bar{x}_{N_1}^{1},\ldots,\bar{x}_1^{j},\ldots, \bar{x}_{i-1}^{j}\bigg)\bigg|^2\bigg]\\
&=\frac{1}{\varepsilon^2 N^2_j}\bar{\mathbb{E}} \bigg[ \sum_{i=1}^{N_j}\bigg(g(\bar{x}_j^j)-\int_{\mathcal{Z}}g(y)\eta_{ji}^N\big(dy |\bar{x}_1^{1},\ldots, \bar{x}_{N_1}^{1},\ldots,\bar{x}_1^{j},\ldots, \bar{x}_{i-1}^{j}\bigg)^2\bigg]\\
&\leq \frac{\|g\|^2_{\infty}}{\varepsilon^2 N_j}.
 \end{align*}
 
 From the last inequality one deduces that, for each $g\in\mathcal{U}_b(\mathcal{Z})$, with $\mathcal{U}_b(\mathcal{Z})$ being the space of bounded and uniformly continuous functions on $\mathcal{Z}$, the sequence $\{\int_{\mathcal{Z}}g(y)d\bar{\mu}^N_j- \int_{\mathcal{Z}}g(y)d\chi_j^N\}_{N\geq 1}$ converges in probability to $0$, and hence in distribution. Moreover, the sequence $\{(\chi_j^N,\bar{\mu}_j^N), N\geq 1\}$ takes values in $\mathcal{M}_1(\mathcal{Z})\times\mathcal{M}_1(\mathcal{Z})$ which is a compact space by the Prokhorov's theorem given that $\mathcal{Z}$ is compact. In addition, $\mathcal{M}_1(\mathcal{Z})\times\mathcal{M}_1(\mathcal{Z})$ endowed with the topology of weak convergence is a metric space. Therefore, any subsequence admits a further subsubsequence that converges weakly to some  $(\chi_j,\bar{\mu}_j)\in\mathcal{M}_1(\mathcal{Z})\times\mathcal{M}_1(\mathcal{Z})$. Thence, by the Skorokhod's representation theorem, there exists some probability space such that, for all $g\in\mathcal{U}_b(\mathcal{Z})$, the following convergences hold almost surely,    
 \begin{align*}
\lim_{N\rightarrow\infty} \int_{\mathcal{Z}}g d\chi_j^N&=\int_{\mathcal{Z}}g d\chi_j,\\
\lim_{N\rightarrow\infty} \int_{\mathcal{Z}}g d\bar{\mu}^N_j&=\int_{\mathcal{Z}}g d\bar{\mu}_j,
 \end{align*}
 and
 \begin{align*}
\lim_{N\rightarrow\infty} \int_{\mathcal{Z}}gd\bar{\mu}^N_j- \int_{\mathcal{Z}}gd\chi_j^N=0.
 \end{align*}
We thus deduce that $\int_{\mathcal{Z}}g d\chi_j=\int_{\mathcal{Z}}g d\bar{\mu}_j$. Then, since the space $\mathcal{U}_b(\mathcal{Z})$ is measure-determining, one deduces that $\chi_j=\bar{\mu}_j$ almost surely. We then conclude that any subsequence of $\{(\chi_j^N,\bar{\mu}_j^N), N\geq 1\}$ contains a further subsubsequence that converge in distribution to $(\bar{\mu}_j,\bar{\mu}_j)$, for some $\bar{\mu}_j\in\mathcal{M}_1(\mathcal{Z})$. Denote by $\{N\in\mathbb{N}\}$ this subsequence. Again by the Skorokhod's representation theorem, there exists some probability space such that $(\chi_j^N,\bar{\mu}_j^N)$ converges almost surely to $(\bar{\mu}_j,\bar{\mu}_j)$ along the subsequence $\{N\in\mathbb{N}\}$. Thence, using the nonnegativy and semi-continuity of the relative entropy function $R(\cdot,\nu)$ on $\mathcal{M}_1(\mathcal{Z})$, the boundedness and continuity of the function $h^N$ on $(\mathcal{M}_1(\mathcal{Z}))^r$ together with its convergence towards $h$, and the convergence $\frac{N_j}{N}\rightarrow \alpha_j$, one obtains
  \begin{align*}
  \liminf_{N\rightarrow\infty}\bigg\{R\bigg(\sum_{j=1}^r\frac{N_j}{N}\chi^N_j\big\|\nu\bigg)+h^N(\bar{\mu}^N)\bigg\}\geq R\bigg(\sum_{j=1}^r\alpha_j\bar{\mu}_j\big\|\nu\bigg)+h(\bar{\mu}).
  \end{align*}
Finally, by Fatou's lemma we obtain
   \begin{align*}
  \liminf_{N\rightarrow\infty} W^N+\varepsilon&\geq   \liminf_{N\rightarrow\infty}\bar{\mathbb{E}}\bigg[ R\bigg( \sum_{j=1}^r \frac{N_j}{N} \chi_j^N(\cdot)\big\| \nu (\cdot)\bigg)+ h^N(\bar{\mu}^N) \bigg]\\
  & \geq \bar{\mathbb{E}}\bigg[ R\bigg(\sum_{j=1}^r\alpha_j\bar{\mu}_j\big\|\nu\bigg)+h^N(\bar{\mu})\bigg]\\
  &\geq \inf_{\gamma\in(\mathcal{M}_1(\mathcal{Z}))^{r}} \bigg\{R\bigg(\sum_{j=1}^r\alpha_j\gamma_j\big\|\nu\bigg)+h(\gamma)\bigg\}.
\end{align*} 
Letting $\varepsilon\downarrow 0$, we conclude that any subsequence of the original sequence $\{W^N,N\in\mathbb{N}\}$ contains a further subsubsequence that satisfies the following lower bound 
 \begin{align}
  \liminf_{N\rightarrow\infty} W^N\geq \inf_{\gamma\in(\mathcal{M}_1(\mathcal{Z}))^{r}} \bigg\{R\bigg(\sum_{j=1}^r\alpha_j\gamma_j\big\|\nu\bigg)+h(\gamma)\bigg\}.
  \label{lower-bound-variat-form} 
 \end{align}
Thence we deduce that the original sequence $\{W^N,N\in\mathbb{N}\}$ satisfies the lower limit in $(\ref{lower-bound-variat-form})$. Combining this with the upper bound $(\ref{upp-bound-variat-form})$ leads to the stated result.   \carre
\\
\\
We are now ready to state the main result of this section.  
 
 \begin{prop}
\label{limit-rel-entr-prop}
There exists a constant $C \in \mathbb{R}$ such that for all $q=(q^{1,c},q^{1,c},\ldots,q^{r,c},q^{r,p})\in\big(\mathcal{M}_1(\mathcal{Z})\big)^{2r}$
\begin{align*}
\lim_{N\rightarrow\infty}\bar{F}^N(q)=F(q),
\end{align*}
where 
\begin{equation}
\begin{split}
\label{Lyap-func2}
F(q)&=\sum_{j=1}^r\bigg[ \frac{\beta}{2} \bigg\{(\alpha_jp_j^c)^2\sum_{z\in\mathcal{Z}}\bigg(\sum_{z'\in\mathcal{Z}}W(z,z')q^{j,c}_{z'}\bigg)q^{j,c}_{z}+2\alpha_j^2p_j^pp_j^c \sum_{z\in\mathcal{Z}}\bigg(\sum_{z'\in\mathcal{Z}}W(z,z')q^{j,p}_{z'}\bigg)q^{j,c}_{z}\\
&\quad+\alpha_jp_j^p\sum_{z\in\mathcal{Z}}\bigg(\sum_{k=1}^r\alpha_kp_k^p\sum_{z'\in\mathcal{Z}}W(z,z')q_{z'}^{k,p}\bigg)q_z^{j,p} \bigg\}\bigg]+\sum_{j=1}^r\bigg( \alpha_jp_j^c \sum_{z\in\mathcal{Z}}q^{j,c}_{z}\log(q^{j,c}_{z})+\alpha_jp_j^p\sum_{z\in\mathcal{Z}}q^{j,p}_{z}\log (q^{j,p}_{z})\bigg)+C.
\end{split}
\end{equation}
\end{prop}  

\proof Fix $q=(q^{1,c},q^{1,p},\ldots,q^{r,c},q^{r,p})\in\big(\mathcal{M}_1(\mathcal{Z})\big)^{2r}$. Then, by $(\ref{stat-distr-gibbs})$ one obtains 
\begin{equation}
\begin{split}
\bar{F}^N(q)&=\frac{1}{N}R\bigg(\otimes^{N_1^c}(q^{1,c})\otimes^{N_1^p}(q^{1,p})\cdots\otimes^{N_r^c}(q^{r,c})\otimes^{N_r^p}(q^{r,p})\|\pi^N({\pmb x})\bigg)\\
            &=\frac{1}{N}\sum_{{\pmb x}\in\mathcal{Z}^N}\prod_{j=1}^r\bigg(\prod_{i\in C_j^c}q^{j,c}_{x_i}\prod_{i\in C_j^p}q^{j,p}_{x_i}\bigg)\log\bigg(\frac{1}{\pi^N({\pmb x})}\prod_{j=1}^r\bigg(\prod\limits_{i\in C_j^c}q^{j,c}_{x_i}\prod\limits_{i\in C_j^p}q^{j,p}_{x_i}\bigg)\bigg)\\
            &=\frac{1}{N}\sum_{{\pmb x}\in\mathcal{Z}^N}\prod_{j=1}^r\bigg(\prod_{i\in C_j^c}q^{j,c}_{x_i}\prod_{i\in C_j^p}q^{j,p}_{x_i}\bigg)\bigg(\sum_{j=1}^r\bigg(\sum_{i\in C_j^c}\log (q^{j,c}_{x_i})+\sum_{i\in C_j^p}\log (q^{j,p}_{x_i})\bigg)\bigg)\\
             &\qquad+\frac{1}{N}\sum_{{\pmb x}\in\mathcal{Z}^N}\prod_{j=1}^r\bigg(\prod_{i\in C_j^c}q^{j,c}_{x_i}\prod_{i\in C_j^p}q^{j,p}_{x_i}\bigg) U_N({\pmb x})+\frac{1}{N}\log Z_N.
\end{split}
\end{equation}

For all $1\leq j\leq r$, let $\{X^{j,c}_i\}_{i\in\mathbb{N}}$ and $\{X^{j,p}_i\}_{i\in\mathbb{N}}$ be sequences of i.i.d. $\mathcal{Z}$-valued random variables with common distributions $q^{j,c}$ and $q^{j,p}$ respectively. We then can write
\begin{align*}
\frac{1}{N}\sum_{{\pmb x}\in\mathcal{Z}^N}\bigg(\prod_{j=1}^r\bigg(\prod_{i\in C_j^c}q^{j,c}_{x_i}\prod_{i\in C_j^p}q^{j,p}_{x_i}\bigg)\bigg)&\bigg(\sum_{j=1}^r\bigg(\sum_{i\in C_j^c} \log(q^{j,c}_{x_i})+\sum_{i\in C_j^p}\log (q^{j,p}_{x_i})\bigg)\bigg)\\
&=\frac{1}{N}\mathbb{E}\bigg(\sum_{j=1}^r\bigg(\sum_{i\in C_j^c} \log(q^{j,c}_{X^{j,c}_i})+\sum_{i\in C_j^p}\log (q^{j,p}_{X^{j,p}_i})\bigg)\bigg)\\
&=\frac{1}{N}\bigg(\sum_{j=1}^r\bigg( N_j^c \mathbb{E}[\log(q^{j,c}_{X^{j,c}_1})]+N_j^p\mathbb{E}[\log (q^{j,p}_{X^{j,p}_1})]\bigg)\bigg)\\
&=\sum_{j=1}^r\bigg( \frac{N_j^c}{N} \sum_{z\in\mathcal{Z}}q^{j,c}_{z}\log(q^{j,c}_{z})+\frac{N_j^p}{N}\sum_{z\in\mathcal{Z}}q^{j,p}_{z}\log (q^{j,p}_{z})\bigg).\\
\end{align*}
Taking the limit in the last display and using Assumption $\ref{ass-prin}$ we find 
\begin{equation}
\begin{split}
\lim_{N\rightarrow\infty}\frac{1}{N}\sum_{{\pmb x}\in\mathcal{Z}^N}\bigg(\prod_{j=1}^r\bigg(\prod_{i\in C_j^c}q^{j,c}_{x_i}\prod_{i\in C_j^p}q^{j,p}_{x_i}\bigg)\bigg)&\bigg(\sum_{j=1}^r\bigg(\sum_{i\in C_j^c} \log(q^{j,c}_{x_i})+\sum_{i\in C_j^p}\log (q^{j,p}_{x_i})\bigg)\bigg)\\
&=\sum_{j=1}^r\bigg( \alpha_j p_j^c \sum_{z\in\mathcal{Z}}q^{j,c}_{z}\log(q^{j,c}_{z})+\alpha_j p_j^p\sum_{z\in\mathcal{Z}}q^{j,p}_{z}\log (q^{j,p}_{z})\bigg).
\label{relat-entr-conv1}
\end{split}
\end{equation}

Moreover, it is easy to see that
\begin{align*}
\frac{1}{N}\sum_{{\pmb x}\in\mathcal{Z}^N}\prod_{j=1}^r\bigg(\prod_{i\in C_j^c}q^{j,c}_{x_i}\prod_{i\in C_j^p}q^{j,p}_{x_i}\bigg)  U_N({\pmb x})=\frac{1}{N}\mathbb{E}[U_N({\pmb X})],
\end{align*} 
where 
 \begin{align*}
U_N({\pmb X})&=\sum_{j=1}^r\bigg[ \frac{\beta}{2N} \bigg\{\sum_{i\in C_j^c}\bigg(\sum_{k\in C^c_j}W(X^{j,c}_i,X^{j,c}_k)+\sum_{k\in C^p_j}W(X^{j,c}_i,X^{j,p}_k)\bigg)\\
&\qquad\qquad+\sum_{i\in C_j^p}\bigg(\sum_{k\in C_j^c} W(X^{j,p}_i,X^{j,c}_k)+\sum_{k\in C_1^p} W(X^{j,p}_i,X^{1,p}_k)+\cdots+\sum_{k\in C_r^p} W(X^{j,p}_i,X^{r,p}_k)\bigg)\bigg\}\bigg],
\end{align*}
is the energy function evaluated at ${\pmb X}=(X^{j,c}_{i_1},X^{j,p}_{i_2},1\leq j\leq r, 1\leq i_1\leq N_j^c, 1\leq i_2\leq N_j^p)$. Thence, relying again on Assumption $\ref{ass-prin}$, one finds
\begin{equation}
\begin{split}
\lim_{N\rightarrow\infty}\frac{1}{N}\mathbb{E}[U_N({\pmb X})]&=\sum_{j=1}^r\bigg[ \frac{\beta}{2} \bigg\{(\alpha_jp_j^c)^2\mathbb{E}[W(X^{j,c}_1,X^{j,c}_2)]+2\alpha_j^2p_j^cp_j^p\mathbb{E}[W(X^{j,c}_1,X^{j,p}_1)]\\
&\quad+\alpha_jp_j^p(\alpha_1p_1^p\mathbb{E}[W(X^{j,p}_1,X^{1,p}_1)]+\cdots+\alpha_jp_j^p\mathbb{E}[W(X^{j,p}_1,X^{j,p}_2)]+\cdots+\alpha_rp_r^p\mathbb{E}[W(X^{j,p}_1,X^{r,p}_1)]\bigg\}\bigg]\\
&=\sum_{j=1}^r\bigg[ \frac{\beta}{2} \bigg\{(\alpha_jp_j^c)^2\sum_{z\in\mathcal{Z}}\bigg(\sum_{z'\in\mathcal{Z}}W(z,z')q^{j,c}_{z'}\bigg)q^{j,c}_{z}+2\alpha_j^2p_j^pp_j^c \sum_{z\in\mathcal{Z}}\bigg(\sum_{z'\in\mathcal{Z}}W(z,z')q^{j,p}_{z'}\bigg)q^{j,c}_{z}\\
&\quad+\alpha_jp_j^p\sum_{z\in\mathcal{Z}}\bigg(\sum_{k=1}^r\alpha_kp_k^p\sum_{z'\in\mathcal{Z}}W(z,z')q_{z'}^{k,p}\bigg)q_z^{j,p} \bigg\}\bigg].\\
\label{relat-entr-conv2}
\end{split}
\end{equation}

Now, in order to evaluate the limit as $N\rightarrow\infty$ of the quantity $\frac{1}{N}\log Z_N$, often referred to as the {\it free energy function} in the literature, let $\{y_i\}_{i\in\mathbb{N}}$ be i.i.d. $\mathcal{Z}$-valued random variables with common distribution $\nu$ given by $\nu_z=\frac{1}{|\mathcal{Z}|}, z  \in\mathcal{Z}$. Then have
\begin{equation}
\begin{split}
Z_N&=\sum_{{\pmb x}\in\mathcal{Z}^N}\exp\{-U_N({\pmb x})\}\\
&=\sum_{{\pmb x}\in\mathcal{Z}^N}\exp\bigg\{- \frac{\beta}{2N} \sum_{j=1}^r\bigg[\sum_{i\in C_j^c}\bigg(\sum_{k\in C^c_j}W(x_i,x_k)+\sum_{k\in C^p_j}W(x_i,x_k)\bigg)\\
&\qquad+\sum_{i\in C_j^p}\bigg(\sum_{k\in C_j^c} W(x_i,x_k)+\sum_{k\in C_1^p} W(x_i,x_k)+\cdots+\sum_{k\in C_r^p} W(x_i,x_k)\bigg)\bigg]\bigg\}\\
&=|\mathcal{Z}|^N\mathbb{E}\bigg[\exp\bigg\{- \frac{\beta}{2N} \sum_{j=1}^r\bigg[\sum_{i\in C_j^c}\bigg(\sum_{k\in C^c_j}W(y_i,y_k)+\sum_{k\in C^p_j}W(y_i,y_k)\bigg)\\
&\qquad+\sum_{i\in C_j^p}\bigg(\sum_{k\in C_j^c} W(y_i,y_k)+\sum_{k\in C_1^p} W(y_i,y_k)+\cdots+\sum_{k\in C_r^p} W(y_i,y_k)\bigg)\bigg]\bigg\} \bigg].
\label{free-energ-estimate}
\end{split}
\end{equation}
Furthermore, recalling $(\ref{loc-emp-meas-gibbs})$, one finds  
\begin{align*}
Z_N&=|\mathcal{Z}|^N\mathbb{E}\bigg[\exp\bigg\{- \frac{\beta}{2N} \sum_{j=1}^r\bigg[N_j^c\sum_{z'\in \mathcal{Z}}\bigg(N_j^c\sum_{z\in \mathcal{Z}}W(z',z)\varrho^{j,c}_z({\pmb y})+N_j^p\sum_{z\in \mathcal{Z}}W(z',z)\varrho^{j,p}_z({\pmb y})\bigg)\varrho^{j,c}_{z'}({\pmb y})\\
&\quad\qquad+N_j^p\sum_{z'\in \mathcal{Z}}\bigg(N_j^c\sum_{z\in \mathcal{Z}} W(z',z)\varrho^{j,c}_z({\pmb y})+N_1^p\sum_{z\in\mathcal{Z}} W(z',z)\varrho_z^{1,p}({\pmb y})+\cdots+N_r^p\sum_{z\in \mathcal{Z}} W(z',z)\varrho_z^{r,p}({\pmb y})\bigg)\varrho^{j,p}_{z'}({\pmb y})\bigg]\bigg\} \bigg].
\end{align*}
Define the sequence of functions $\{\Phi^N,N\in\mathbb{N}\}$ mapping $\big(\mathcal{M}_1(\mathcal{Z})\big)^{2r}$ into $\mathbb{R}$ by
\begin{align*}
\Phi^N(\eta)   &= \frac{\beta}{2N^2} \sum_{j=1}^r\bigg[N_j^c\sum_{z'\in \mathcal{Z}}\bigg(N_j^c\sum_{z\in \mathcal{Z}}W(z',z)\eta^{j,c}_z+N_j^p\sum_{z\in \mathcal{Z}}W(z',z)\eta^{j,p}_z\bigg)\eta^{j,c}_{z'}\\
&\quad+N_j^p\sum_{z'\in \mathcal{Z}}\bigg(N_j^c\sum_{z\in \mathcal{Z}} W(z',z)\eta^{j,c}_z+N_1^p\sum_{z\in\mathcal{Z}} W(z',z)\eta_z^{1,p}+\cdots+N_r^p\sum_{z\in \mathcal{Z}} W(z',z)\eta_z^{r,p}\bigg)\eta^{j,p}_{z'}\bigg],
\end{align*}
 for any $\eta=(\eta^{1,c},\eta^{1,p},\ldots,\eta^{r,c},\eta^{r,p})\in\big(\mathcal{M}_1(\mathcal{Z})\big)^{2r}$. Notice that by Assumption \ref{ass-prin}, it is easy to see that the functions $\Phi^N$ are bounded continuous for all $N$ and converge towards the function $\Phi$ given by
\begin{align*}
\Phi(\eta)   &= \frac{\beta}{2} \sum_{j=1}^r\bigg[p_j^c\sum_{z'\in \mathcal{Z}}\bigg(p_j^c\sum_{z\in \mathcal{Z}}W(z',z)\eta^{j,c}_z+p_j^p\sum_{z\in \mathcal{Z}}W(z',z)\eta^{j,p}_z\bigg)\eta^{j,c}_{z'}\\
&\quad+p_j^p\sum_{z'\in \mathcal{Z}}\bigg(p_j^c\sum_{z\in \mathcal{Z}} W(z',z)\eta^{j,c}_z+p_1^p\sum_{z\in\mathcal{Z}} W(z',z)\eta_z^{1,p}+\cdots+p_r^p\sum_{z\in \mathcal{Z}} W(z',z)\eta_z^{r,p}\bigg)\eta^{j,p}_{z'}\bigg].
\end{align*}
Hence, using the Laplace principle given in Proposition \ref{Sanov-vect-meas-prop}, the following convergence holds true
\begin{align}
\lim_{N\rightarrow\infty}\frac{1}{N}\log Z_N=-\inf_{\gamma\in(\mathcal{M}_1(\mathcal{Z}))^{r}} \bigg\{R\bigg(\sum_{j=1}^r(p^c_j\gamma^c_j(\cdot)+p^p_j\gamma^p_j(\cdot))\big\|\nu(\cdot)\bigg)+\Phi(\gamma)\bigg\}+\log\mathcal{|Z|}=C,
\label{relat-entr-conv3}
\end{align}
where $C$ is a finite constant that does not depends on the measure $q$. Combining $(\ref{relat-entr-conv1}),(\ref{relat-entr-conv2})$, and $(\ref{relat-entr-conv3})$ gives $(\ref{Lyap-func2})$. The proof is complete. 
\carre

\subsection{Fixed points of the McKean-Vlasov system}
One of the important features of the limiting function $F(q)$ given by $(\ref{Lyap-func2})$ is that it characterizes the critical points of the McKean-Vlasov system in $(\ref{McKean-Vlas-syst})$. To prove this fact, let us introduce some additional notations.  Define the hyperplane
\begin{align*}
\mathcal{H}^1=\bigg\{v\in(\mathbb{R}^K)^{2r}:\sum_{i=1}^Kv_i^{j,\iota}=1,\mbox{ for all $1\leq j\leq r$, $\iota\in\{c,p\}$}\bigg\},
\end{align*}
and the corresponding shifted version
\begin{align*}
\mathcal{H}^0=\bigg\{v\in(\mathbb{R}^K)^{2r}:\sum_{i=1}^Kv_i^{j,\iota}=0,\mbox{ for all $1\leq j\leq r$, $\iota\in\{c,p\}$}\bigg\}.
\end{align*}
For any $v\in\mathcal{H}^0$, the directional derivative of the function $F(q)$ with respect to $v$ is given by
\begin{align*}
\frac{\partial}{\partial v}F(q)=\sum_{j=1}^r\sum_{\iota\in\{c,p\}}\bigg[\sum_{x\in\mathcal{Z}}\bigg(\frac{\partial F(q)}{\partial q_x^{j,\iota}}\bigg)v_x^{j,\iota}\bigg],
\end{align*}
where, for $x\in\mathcal{Z}$ and $1\leq j\leq r$,
\begin{align*}
\frac{\partial}{\partial q_x^{j,c}}F(q)&=\beta\bigg[ (\alpha_jp_j^c)^2 \sum_{z\in\mathcal{Z}}W(x,z)q^{j,c}_{z}+\alpha_j^2p_j^pp_j^c \sum_{z\in\mathcal{Z}}W(x,z)q^{j,p}_{z}\bigg]+ \alpha_jp_j^c (\log(q^{j,c}_{x})+1),
\end{align*}
and
\begin{align*}
\frac{\partial}{\partial q_x^{j,p}}F(q)&=\beta\bigg[\alpha_j^2p_j^pp_j^c \sum_{z\in\mathcal{Z}}W(z,x)q^{j,c}_{z}+\alpha_jp_j^p\sum_{k=1}^r\alpha_kp_k^p\sum_{z\in\mathcal{Z}}W(x,z)q_{z}^{k,p}\bigg]+ \alpha_jp_j^p (\log (q^{j,p}_{x})+1).
\end{align*}

\begin{lem}
Let the rate functions $\lambda^{j,c}_{z,z'}$ and $\lambda^{j,p}_{z,z'}$ be given by $(\ref{lambd-gibbs})$, and consider the corresponding McKean-Vlasov system in $(\ref{McKean-Vlas-syst})$. Then, a given $q=(q^{1,c},q^{1,p},\ldots,q^{r,c},q^{r,p})\in(\mathcal{M}_1(\mathcal{Z}))^{2r}$ is a fixed point of $(\ref{McKean-Vlas-syst})$ if and only if $q\in \mathcal{S}^{\circ^{2r}}$  and $\frac{d}{dv}F(q)=0$ for all $v\in\mathcal{H}^0$. 
\label{crit-point-lem}
\end{lem} 
 
 \proof Define the  vector $\bar{v}=(e_x-e_y,\ldots,e_x-e_y)\in\mathcal{H}^0$, where $e_x$ is the unit vector of $\mathbb{R}^K$ in the $x$-direction. Note that given the structure of the shifted space $\mathcal{H}^{\circ}$, it is enough to prove the result for $\bar{v}$. First, the directional derivative of $F(q)$ with respect to $\bar{v}$ is given by 
\begin{align*}
\frac{\partial}{\partial \bar{v}}F(q)=\sum_{j=1}^r&\bigg[\alpha_jp_j^c \bigg( \beta \bigg\{\alpha_jp_j^c \sum_{z\in\mathcal{Z}}\big(W(x,z)-W(y,z)\big)q^{j,c}_{z}+\alpha_jp_j^p \sum_{z\in\mathcal{Z}}\big(W(x,z)-W(y,z)\big)q^{j,p}_{z}\bigg\}\\
&+  (\log(q^{j,c}_{x})-\log(q^{j,c}_{y}))\bigg)\\
&+\alpha_jp_j^p\bigg(\beta\bigg\{\alpha_jp_j^c \sum_{z\in\mathcal{Z}}\big(W(z,x)-W(z,y)\big)q^{j,c}_{z}+\sum_{k=1}^r\alpha_kp_k^p\sum_{z\in\mathcal{Z}}\big(W(z,x)-W(z,y)\big)q_{z}^{k,p}\bigg\}\\
&+ (\log (q^{j,p}_{x})-\log (q^{j,p}_{y}))\bigg)
\bigg].
\end{align*} 
Moreover by $(\ref{stat-dist-fix-c})$ and $(\ref{stat-dist-fix-p})$ one can rewrite the last equality as
\begin{align}
\frac{\partial}{\partial \bar{v}}F(q)=\sum_{j=1}^r&\bigg[\alpha_jp_j^c \bigg( \log\bigg(\frac{q^{j,c}_{x}}{q^{j,c}_{y}}\bigg) -\log\bigg(\frac{\pi^{j,c}(q)_x}{\pi^{j,c}(q)_y}\bigg)\bigg)+\alpha_jp_j^p\bigg(\log\bigg(\frac{q^{j,p}_{x}}{q^{j,p}_{y}}\bigg) -\log\bigg(\frac{\pi^{j,p}(q)_x}{\pi^{j,p}(q)_y}\bigg)\bigg)
\bigg].
\label{direc-deriv-unit}
\end{align}
 
Let $q=(q^{1,c},q^{1,p},\ldots,q^{r,c},q^{r,p})\in\big(\mathcal{M}_1(\mathcal{Z})\big)^{2r}$ be a fixed point of the McKean-Vlasov system in $(\ref{McKean-Vlas-syst})$ corresponding to the rate functions $\lambda^{j,c}_{z,z'}$ and $\lambda^{j,p}_{z,z'}$ defined in $(\ref{lambd-gibbs})$. Thus, by the definition, 
\begin{equation}
\begin{split}
\left\{ \begin{array}{lcl}q^{j,c}A^{j,c}_q=0, & &  \\
q^{j,p}A^{j,p}_q=0,  & & \\
 1\leq j\leq r. & &
 \end{array}\right.
\end{split}
\end{equation} 
Moreover, from $(\ref{stat-dist-fix-c})$ and $(\ref{stat-dist-fix-p})$ and $\pi^{j,c}(q),\pi^{j,p}(q)$ are the stationary distributions of irreducible and recurrent continuous-time Markov chains, therefore, 
\begin{equation}
\begin{split}
\left\{ \begin{array}{lcl}\pi^{j,c}_qA^{j,c}_q=0, & &  \\
\pi^{j,p}_qA^{j,p}_q=0,  & & \\
 1\leq j\leq r. & &
 \end{array}\right.
\end{split}
\end{equation} 
But, since the corresponding Markov chains are irreducible and positive recurrent, there exists a unique solution to the associated balance equations, and thus $\pi^{j,\iota}(q)=q^{j,\iota}$ for all $\iota\in\{c,p\}$ and $1\leq j\leq r$. Furthermore, since the stationary distribution satisfies $\pi^{j,\iota}(q)_x>0$, we have $q^{j,\iota}_x>0$ for all $x\in\mathcal{Z}$, $\iota\in\{c,p\}$ and $1\leq j\leq r$. Thus $q^{j,\iota}\in\mathcal{S}^{\circ}$ for all $\iota\in\{c,p\}$ and $1\leq j\leq r$ and thence $q\in\mathcal{S}^{\circ^{2r}}$. In addition, using $(\ref{direc-deriv-unit})$ we find that $\frac{\partial}{\partial \bar{v}}F(q)=0$. 

Conversely, suppose that $\frac{\partial}{\partial \bar{v}}F(q)=0$. Then from $(\ref{direc-deriv-unit})$ we straightforwardly obtain that
\begin{align*}
\bigg(\frac{q^{j,\iota}_{x}}{q^{j,\iota}_{y}}\bigg)=\bigg(\frac{\pi^{j,\iota}(q)_x}{\pi^{j,\iota}(q)_y}\bigg),\mbox{
for all $x,y\in\mathcal{Z}$, $\iota\in\{c,p\}$ and $1\leq j\leq r$.}
\end{align*}
 Thence $q=\pi(q)$, and thus $q$ is a fixed point of the McKean-Vlasov system. This concludes the proof. \carre 
\\
\\
The previous result allows us to identify the equilibrium points of the McKean-Vlasov system in $(\ref{McKean-Vlas-syst}) $ by the critical points of the limit function $ F $. Also, note that the dynamic system can contain multiple $\omega $ -limit sets as shown in the following example.
\begin{ex}
\label{exe-multip-omega}
Suppose $\mathcal{Z}=\{1,2\}$, $W(1,1)=W(2,2)=0$ and $W(1,2)=W(2,1)=1$. Moreover, suppose that $r=2$, $\alpha_1=\alpha_2=p_1^c=p_1^p=p_2^c=p_2^p=\frac{1}{2}$. Therefore, by $(\ref{Lyap-func2})$, $F(q)=f(q^{1,c}_{1},q^{1,p}_{1},q^{2,c}_{1},q^{2,p}_{1})$, where 
\begin{align*}
f(q^{1,c}_{1},q^{1,p}_{1},q^{2,c}_{1},q^{2,p}_{1})&=\sum_{j=1}^2\bigg[\frac{\beta}{2} \bigg\{2\bigg(\frac{1}{4}\bigg)^2 q^{j,c}_{1}(1-q^{j,c}_{1})+2\bigg(\frac{1}{4}\bigg)^2 \big(q^{j,c}_{1}(1-q^{j,p}_{1})+q^{j,p}_{1}(1-q^{j,c}_{1})\big)\\
&\qquad\qquad+\frac{1}{4}\bigg(q^{j,p}_1\sum_{k=1}^2\frac{1}{4}(1-q_1^{k,p})+(1-q^{j,p}_1)\sum_{k=1}^2\frac{1}{4} q_1^{k,p} \bigg)\bigg\}\bigg]\\
&\quad+\sum_{j=1}^2\bigg[ \frac{1}{4} \big(q^{j,c}_{1}\log(q^{j,c}_{1})+(1-q^{j,c}_{1})\log(1-q^{j,c}_{1})\big)+\frac{1}{4} \big(q^{j,p}_{1}\log(q^{j,p}_{1})+(1-q^{j,p}_{1})\log(1-q^{j,p}_{1})\big)\bigg].
\end{align*}
Hence, the critical points of $f(q^{1,c}_{1},q^{1,p}_{1},q^{2,c}_{1},q^{2,p}_{1})$ on $[0,1]^4$ correspond to the critical points of $F$ on $(\mathcal{M}_1(\mathcal{Z}))^{4}$ which, by Lemma \ref{crit-point-lem}, correspond to the equilibria  of the McKean-Vlasov system $(\ref{McKean-Vlas-syst})$. Therefore, one has to solve the following system of equations: 
\begin{equation}
\begin{split}
\left\{\begin{array}{c}
\frac{\partial}{\partial q^{1,c}_1}f(q^{1,c}_{1},q^{1,p}_{1},q^{2,c}_{1},q^{2,p}_{1})=0,
\\
\frac{\partial}{\partial q^{1,p}_1}f(q^{1,c}_{1},q^{1,p}_{1},q^{2,c}_{1},q^{2,p}_{1})=0,\\
\frac{\partial}{\partial q^{2,c}_1}f(q^{1,c}_{1},q^{1,p}_{1},q^{2,c}_{1},q^{2,p}_{1})=0,\\
\frac{\partial}{\partial q^{2,p}_1}f(q^{1,c}_{1},q^{1,p}_{1},q^{2,c}_{1},q^{2,p}_{1})=0.\\
\end{array}\right.
\label{parti-syst}
\end{split}
\end{equation}
By straightforward calculations, the partial derivatives of $f$ are given by
\begin{align*}
\frac{\partial}{\partial q^{j,c}_1}f(q^{1,c}_{1},q^{1,p}_{1},q^{2,c}_{1},q^{2,p}_{1})&=\beta \bigg(\frac{1}{4}\bigg)^2 \bigg(1-2q^{j,c}_{1}\bigg)+ \beta \bigg(\frac{1}{4}\bigg)^2  \bigg(1-2q_1^{j,p}\bigg)+\frac{1}{4} \bigg(\log(q^{j,c}_{1})-\log(1-q^{j,c}_{1})\bigg),
\\
\frac{\partial}{\partial q^{j,p}_1}f(q^{1,c}_{1},q^{1,p}_{1},q^{2,c}_{1},q^{2,p}_{1})&=\beta \bigg(\frac{1}{4}\bigg)^2 \bigg(1-2q^{j,c}_{1}\bigg)+ \beta \frac{1}{4} \bigg(\sum_{k=1}^2\frac{1}{4}(1-2q_1^{k,p})\bigg)+\bigg(\frac{1}{4}\bigg)^2  \bigg(\log(q^{j,p}_{1})-\log(1-q^{j,p}_{1})\bigg),
\end{align*}
 for all $1\leq j\leq 2$. Fixing $\beta=4$, the numerical resolution of the system of equations in $(\ref{parti-syst})$ identifies the following three solutions: 
\begin{align*}
(q^{1,c}_{1},q^{1,p}_{1},q^{2,c}_{1},q^{2,p}_{1})&=(\frac{1}{2},\frac{1}{2},\frac{1}{2},\frac{1}{2}),
\\
(q^{1,c}_{1},q^{1,p}_{1},q^{2,c}_{1},q^{2,p}_{1})&=(0.8039,0.9015,0.8039,0.9015),\\
(q^{1,c}_{1},q^{1,p}_{1},q^{2,c}_{1},q^{2,p}_{1})&=(0.1961,0.0985,0.1961,0.0985).\\
\end{align*}
Consequently, these solutions correspond to the following critical points of the function $F(q)$:
 \begin{align*}
(q^{1,c}_{1},q^{1,c}_{2},q^{1,p}_{1},q^{1,p}_{2},q^{2,c}_{1},q^{2,c}_{2},q^{2,p}_{1},q^{2,p}_{2})&=(\frac{1}{2},\frac{1}{2},\frac{1}{2},\frac{1}{2},\frac{1}{2},\frac{1}{2},\frac{1}{2},\frac{1}{2}),
\\
(q^{1,c}_{1},q^{1,c}_{2},q^{1,p}_{1},q^{1,p}_{2},q^{2,c}_{1},q^{2,c}_{2},q^{2,p}_{1},q^{2,p}_{2})&=(0.8039,0.1961,0.9015,0.0985,0.8039,0.1961,0.9015,0.0985),\\
(q^{1,c}_{1},q^{1,c}_{2},q^{1,p}_{1},q^{1,p}_{2},q^{2,c}_{1},q^{2,c}_{2},q^{2,p}_{1},q^{2,p}_{2})&=(0.1961,0.8039,0.0985,0.9015,0.1961,0.8039,0.0985,0.9015),
\end{align*}
which, following Lemma \ref{crit-point-lem}, correspond to the critical points of the McKean-Vlasov system in $(\ref{McKean-Vlas-syst})$. 

\end{ex}

\subsection{Descent property and Lyapunov function}
In this final section we show that, under positive definitness assumption, the limiting function $F(q)$ in $(\ref{Lyap-func2})$ is indeed a local Lyapunov function for the McKean-Vlasov system in $(\ref{McKean-Vlas-syst})$. The next result proves that $F(q)$ satisfies a descent property.

\begin{prop}
\label{descent-prop}
Let $\mu(t)$ be the solution to the McKean-Vlasov system in $(\ref{McKean-Vlas-syst})$ corresponding to the rate functions $\lambda_{z,z'}^{j,c}$ and $\lambda_{z,z'}^{j,p}$ defined in $(\ref{lambd-gibbs})$, and starting at some $\mu(0)\in(\mathcal{M}_1(\mathcal{Z}))^{2r}$. Then, for all $t> 0$, we have
\begin{align*}
\frac{d}{dt}F(\mu(t))&=\frac{d}{dt}\sum_{j=1}^r\bigg(\alpha_jp_j^c R\big(\mu^{j,c}\|\pi^{j,c}(\nu^{j,c})\big)+\alpha_jp_j^p R\big(\mu^{j,p}\|\pi^{j,p}(\nu^{j,p})\big)\bigg)\bigg|_{\nu=\mu(t)}\leq 0. 
\end{align*}
  
  Moreover, $\frac{d}{dt}F(\mu(t))=0$ if and only if $\mu^{j,\iota}(t)=\pi^{j,\iota}(\mu^{j,\iota}(t))$ for all $\iota\in\{c,p\}$ and $1\leq j\leq r$, where $\pi^{j,c}(\cdot)$ and $\pi^{j,p}(\cdot)$ are defined in $(\ref{stat-dist-fix-c})$ and $(\ref{stat-dist-fix-p})$, respectively.
  \label{descent-prop-lyap}
\end{prop}

\proof First, recall that any probability flow $q(t)$ on $\mathcal{Z}$ satisfies $\frac{d}{dt}\sum_{z\in\mathcal{Z}}q_z(t)=\sum_{z\in\mathcal{Z}}\frac{d}{dt}q_z(t)=0$ for all $t\geq 0$. Using this together with $(\ref{Lyap-func2})$ one obtains 
\begin{equation}
\begin{split}
\frac{d}{dt}F(\mu(t))&=\sum_{j=1}^r\bigg[ \frac{\beta}{2} \bigg\{(\alpha_jp_j^c)^2\sum_{z\in\mathcal{Z}}\bigg(\bigg(\sum_{z'\in\mathcal{Z}}W(z,z')\mu^{j,c}_{z'}(t)\bigg)\frac{d\mu^{j,c}_{z}}{dt}(t)+ \bigg(\sum_{z'\in\mathcal{Z}}W(z,z')\frac{d\mu^{j,c}_{z'}}{dt}(t)\bigg)\mu^{j,c}_{z}(t)\bigg)\\
&\quad+2\alpha_j^2p_j^pp_j^c \sum_{z\in\mathcal{Z}}\bigg(\bigg(\sum_{z'\in\mathcal{Z}}W(z,z')\mu^{j,p}_{z'}(t)\bigg)\frac{d\mu^{j,c}_{z}}{dt}(t)+\bigg(\sum_{z'\in\mathcal{Z}}W(z,z')\frac{d\mu^{j,p}_{z'}}{dt}(t)\bigg)\mu^{j,c}_{z}(t)\bigg)\\
&\quad+\alpha_jp_j^p\sum_{z\in\mathcal{Z}}\bigg(\bigg(\sum_{k=1}^r\alpha_kp_k^p\sum_{z'\in\mathcal{Z}}W(z,z')\frac{d\mu_{z'}^{k,p}}{dt}(t)\bigg)\mu_z^{j,p}(t)+\bigg(\sum_{k=1}^r\alpha_kp_k^p\sum_{z'\in\mathcal{Z}}W(z,z')\mu_{z'}^{k,p}(t)\bigg)\frac{d\mu_z^{j,p}}{dt}(t)\bigg) \bigg\}\bigg]\\
&\quad+\sum_{j=1}^r\bigg( \alpha_jp_j^c \sum_{z\in\mathcal{Z}}\frac{d\mu^{j,c}_{z}}{dt}(t)\log(\mu^{j,c}_{z}(t))+\alpha_jp_j^p\sum_{z\in\mathcal{Z}}\frac{d\mu^{j,p}_{z}}{dt}(t)\log (\mu^{j,p}_{z}(t))\bigg).
\end{split}
\end{equation}
In addition, by $(\ref{stat-dist-fix-c})$ and $(\ref{relat-entro-def})$ one further obtains, for any $\nu=(\nu^{1,c},\nu^{1,p},\ldots,\nu^{r,c},\nu^{r,p})\in(\mathcal{M}_1(\mathcal{Z}))^{2r}$,
\begin{align*}
\frac{d}{dt}R\bigg(\mu^{j,c}(t)\big\|\pi^{j,c}(\nu^{j,c})\bigg)=\sum_{z\in\mathcal{Z}}\frac{d\mu_z^{j,c}}{dt}(t)\log\mu_z^{j,c}(t)+\beta\sum_{z\in\mathcal{Z}}\bigg(& \alpha_jp_j^c\sum_{z'\in\mathcal{Z}}W(z,z')\nu^{j,c}_{z'}\\
&+\alpha_jp_j^p\sum_{z'\in\mathcal{Z}}W(z,z')\nu^{j,p}_{z'}\bigg)\frac{d\mu^{j,c}_{z}}{dt}(t).
\end{align*}
Moreover, using $(\ref{stat-dist-fix-p})$ and the definition of the relative entropy one finds 
\begin{align*}
\frac{d}{dt}R\bigg(\mu^{j,p}(t)\big\|\pi^{j,p}(\nu^{j,p})\bigg)=\sum_{z\in\mathcal{Z}}\frac{d\mu_z^{j,p}}{dt}(t)\log\mu_z^{j,p}(t)+\beta\sum_{z\in\mathcal{Z}}\bigg(& \alpha_jp_j^c\sum_{z'\in\mathcal{Z}}W(z,z')\nu^{j,c}_{z'}\\
&+\sum_{l=1}^r\alpha_lp_l^p\sum_{z'\in\mathcal{Z}}W(z,z')\nu^{l,p}_{z'}\bigg)\frac{d\mu^{j,p}_{z}}{dt}(t).
\end{align*}
Thence, since $\nu$ is arbitrary, one observes that, for all $t\geq 0$,  
\begin{align*}
\frac{d}{dt}F(\mu(t))&=\frac{d}{dt}\sum_{j=1}^r\bigg(\alpha_jp_j^c R\big(\mu^{j,c}(t)\|\pi^{j,c}(\nu^{j,c})\big)+\alpha_jp_j^p R\big(\mu^{j,p}(t)\|\pi^{j,p}(\nu^{j,p})\big)\bigg)\bigg|_{\nu=\mu(t)}.
\end{align*}

Moreover, recall that, since $\nu$ is fixed, $\pi^{j,c}(\nu^{j,c})$ and $\pi^{j,p}(\nu^{j,p})$ are the stationary distributions of ergodic Markov processes generated by the rate functions $(\ref{lambd-gibbs})$. The linear Kolmogorov forward equations associated with these Markov processes are given by 
\begin{equation}
\begin{split}
\left\{ \begin{array}{lcl}\dot{\eta}_j^c(t)=\eta_j^c(t)A^{j,c}_{\nu}, & &  \\
 \dot{\eta}_j^p(t)=\eta_j^p(t)A^{j,p}_{\nu},  & & \\
 1\leq j\leq r, t\geq 0.  & &
 \end{array}\right.
\end{split}
\end{equation}
Fix $t\geq 0$. Since $\nu$ is arbitrary, then one can take $\nu=\mu(t)$ in the last equation. Therefore, at time $t$, both $\mu(t)=(\mu^{1,c}(t),\mu^{1,p}(t),\ldots,\mu^{r,c}(t),\mu^{r,p}(t))$ and $\eta(t)=(\eta^{1,c}(t),\eta^{1,p}(t),\ldots,\eta^{r,c}(t),\eta^{r,p}(t))$ solve a Kolmogorov forward equation with the same rate matrix $A_{\mu(t)}$. Moreover, recall that the relative entropy function  has descent property along the solution to the linear Kolmogorov forward equation. In particular, the proof of this relies on the fact that two solutions to the linear equation with different initial conditions satisfy the forward equation with the same fixed rate matrix; See, e.g., the proof in \cite[Lem 3.1]{Budh+Du(a)2015}. Therefore, since this fact is also satisfied here the descent property of the relative entropy function still holds true and thus    
 \begin{align*}
\frac{d}{dt}F(\mu(t))&=\frac{d}{dt}\sum_{j=1}^r\bigg(\alpha_jp_j^c R\big(\mu^{j,c}(t)\|\pi^{j,c}(\nu^{j,c})\big)+\alpha_jp_j^p R\big(\mu^{j,p}(t)\|\pi^{j,p}(\nu^{j,p})\big)\bigg)\bigg|_{\nu=\mu(t)}\leq 0.
\end{align*}
 
 Using again \cite[Lem 3.1.]{Budh+Du(a)2015} gives that $\frac{d}{dt}F(\mu(t))=0$ if and only if $\mu^{j,\iota}(t)=\pi^{j,\iota}(\mu^{j,\iota}(t))$, for all $\iota\in\{c,p\}$ and $1\leq j\leq r$. \carre
 \\
 \\
The descent property established in Proposition \ref{descent-prop} together with the observation that $F(q)$ is $\mathcal{C}^1$ and uniformly continuous shows that, if $F$ is positive definite in the neighborhood $\Gamma$ of a fixed point $\phi$ of the McKean-Vlasov system in $(\ref{McKean-Vlas-syst})$, then $F(q)$ is a local Lyapunov function in the sense of Definition $\ref{lyapu-def}$. In such a case, the corresponding fixed point $\phi$ is locally stable by Proposition $\ref{stab-lyap}$.

\begin{rem}
\begin{enumerate}
\item If the McKean-Vlasov system in $(\ref{McKean-Vlas-syst})$ contains multiple fixed points, then under positive definiteness assumptions, the function $F$ serves as a Lyapunov function for all of these points and thus allows studying their local stability. Note that if the McKean-Vlasov system in $(\ref{McKean-Vlas-syst})$ contains multiple (local) stable equilibria, one can investigate the metastability of the corresponding $N$-particles system. It amounts to studying the transitions of the empirical vector $\mu^N$ between the different attractor states as time becomes large. Note that these transitions happen even though the existence of a unique invariant measure for $\mu^N$. It is then of interest to estimate quantities such as the mean time spent by the process near a stable point, the probability of reaching a given stable point before reaching another one, or also the probability of transiting between a collection of $\omega$-limit sets in a particular order, and so on. Notice that the metastability analysis was conducted in \cite{Daw+Sid+zha2021} for the general finite-state mean-field systems on block graphs under the classical conditions of Friedlin and Wentzell \cite{Freid+Wentz2012}. 
\item In the general case studied in \cite{Daw+Sid+zha2020,Daw+Sid+zha2021}, the explicit form of the invariant measure is in general not available. Therefore, one cannot evaluate the $N\rightarrow\infty$ limit of the function $\bar{F}^N(q)$ defined in $(\ref{Lyap-func})$. An alternative approach consists in evaluating the limit of 
\begin{align*}
F^N_t(q)=\frac{1}{N}R\bigg(\otimes^{N_1^c}q^{1,c}\otimes^{N_1^p}q^{1,p}\cdots\otimes^{N_r^c}q^{r,c}\otimes^{N_r^p}q^{r,p}\bigg\|  {\pmb p}^N (t)\bigg),
\end{align*}
as $N\rightarrow\infty$ and $t\rightarrow\infty$ where ${\pmb p}^N(t)$ is a multi-exchangeable probability distribution of ${\pmb x}^N(t)=\big(x_n(t),x_m(t),n\in C_j^c, m\in C_j^p,1\leq j\leq r\big)$ at time $t$. Then, providing that this limit takes a useful form, one can investigate its Lyapunov properties. Note that this approach was shown to work in  \cite{Budh+Du(a)2015,Budh+Du(b)2015} for some family of finite-state mean-field models on complete graphs. The main condition is the existence of large deviations principle for the empirical measure. It is then of interest to apply similar strategy to our multi-class setting. One need to rely on the large deviations principles for the empirical vector process, which was proved to hold  in \cite{Daw+Sid+zha2021}. However this goes beyond the scope of this paper. 
\end{enumerate} 
\end{rem}

  \section*{Acknowledgment}
This research was supported by the Natural Sciences and Engineering Research Council of Canada Discovery Grants and by Carleton University.

\bibliographystyle{livre} 
\bibliography{biblio}
\end{document}